\documentclass[11pt]{article}

\usepackage{amsmath,amsthm}

\usepackage{amssymb,latexsym}

\usepackage{enumerate}

\newcommand{\dyw}{\mbox{\rm div}}

\newtheorem{tw}{Theorem}[subsection]

\newtheorem{lm}[tw]{Lemma}
\newtheorem{uw}[tw]{Remark}
\newtheorem{wn}[tw]{Corollary}
\newtheorem{stw}[tw]{Proposition}
\newenvironment{dow}{\it Proof.\rm}{\hfill $\Box$}

\newcommand{\nsubsection}{\setcounter{equation}{0}\subsection}

\begin{document}

\title {Strong solutions of semilinear parabolic equations with measure
data and generalized backward stochastic differential
equations\footnote{Research supported by the Polish Minister of
Science and Higher Education under Grant N N201 372 436. }}
\author {Tomasz Klimsiak \smallskip\\
{\small Faculty of Mathematics and Computer Science,
Nicolaus Copernicus University} \\
{\small  Chopina 12/18, 87--100 Toru\'n, Poland}\\
{\small e-mail: tomas@mat.uni.torun.pl}}
\date{}
\maketitle \noindent{\small{\bf Abstract:} We prove that under
natural assumptions on the data  strong solutions in Sobolev
spaces of semilinear parabolic equations in divergence form
involving  measure on the right-hand side may be represented by
solutions of some generalized backward stochastic differential
equations. As an application we provide stochastic representation
of strong  solutions of the obstacle problem be means of solutions
of some reflected backward stochastic differential equations. To
prove the latter result we use a stochastic homographic
approximation for solutions of the reflected backward equation.
The approximation may be viewed as a stochastic analogue of the
homographic approximation for solutions to the obstacle problem.}
\medskip\\
{\small{\bf Keywords:} Semilinear parabolic equation, Measure
data, Obstacle problem, Strong solution, Generalized BSDE,
Reflected BSDE, Homographic approximation.}
\medskip\\
{\small{\bf 2000 Mathematics Subject Classification:} Primary
35K85, 60H30, Secondary 60H35}

\nsubsection{Introduction}

Let $\mu$ be a Radon measure on $Q_T\equiv[0,T]\times\mathbb{R}^d$
and let $\varphi:\mathbb{R}^d\rightarrow\mathbb{R}$,
$f:Q_T\times\mathbb{R} \times\mathbb{R}^d\rightarrow\mathbb{R}$,
$g:Q_T\times\mathbb{R} \rightarrow\mathbb{R}$ be measurable
functions. In the paper we consider strong solutions in Sobolev
spaces of the Cauchy problem
\begin{equation}
\label{eqm.0} \frac{\partial u}{\partial
t}+L_{t}u=-f_{u}-g(u)\mu,\quad u(T)=\varphi.
\end{equation}
Here
\begin{equation}\label{eq1.3}
L_{t}=\frac{1}{2}\sum_{i,j=1}^d\frac{\partial}{\partial
x_{i}}(a^{ij}\frac{\partial}{\partial x_{j}})+\sum_{i=1}^d
b^{i}\frac{\partial}{\partial x_{i}}
\end{equation}
is an operator with measurable coefficients $a:Q_T
\rightarrow\mathbb{R}^{d}\otimes\mathbb{R}^{d}$,
$b:Q_T\rightarrow\mathbb{R}^{d}$  such that
\begin{equation}\label{eq1.1}
\lambda|\xi|^2\leq\sum^d_{i,j=1}a^{ij}(t,x)\xi_i\xi_j
\le\Lambda|\xi|^2,\quad a^{ij}=a^{ji}, \quad|b^{i}(t,x)|\le
\Lambda,\quad\xi\in\mathbb{R}^d
\end{equation}
for some $0<\lambda\le\Lambda$, and
$f_u(t,x)=f(t,x,u(t,x),\sigma\nabla u(t,x))$ with $\sigma$ such
that $\sigma\sigma^*=a$, $g(u)(t,x)=g(t,x,u(t,x))$, $(t,x)\in
Q_T$.

Nonlinear elliptic equations in divergence form with measure data
on the right-hand side are considered in \cite{Murat}. Our
interest in general parabolic equations  of the form (\ref{eqm.0})
with nonlinear $g$ comes from that fact that as we shall see in
Section \ref{sec5} they arise naturally when considering the
obstacle problem (parabolic variational inequalities). Let us
mention also that equations of the form (\ref{eqm.0}) include the
so-called Schr\"odinger equations with measure data, that is
parabolic equations of the form (\ref{eqm.0}) with $g(t,x,y)=y$
considered for example in \cite{Getoor}. Abstract  parabolic
evolution equations involving measures which depend nonlinearly on
the solution are considered in \cite{Amann}.

Let $\mathbb{X}=\{(X,P_{s,x});(s,x)\in[0,T)\times{\mathbb{R}}^d\}$
be a Markov family corresponding to the operator $L_{t}$ (see
\cite{Ro1,Stroock.DIFF}). Our main result concerning (\ref{eqm.0})
says that if $\mu$ belongs to the weighted  Sobolev space
${\mathbb{L}}_{2}(0,T;H^{-1}_\varrho)$ then under
natural conditions on $\varphi,f,g$ there exists a minimal strong
solution of (\ref{eqm.0}) and the pair
\begin{equation}
\label{eq1.04} (Y^{s,x}_t,Z^{s,x}_t)=(u(t,X_t),\sigma\nabla
u(t,X_t)),\quad t\in[s,T]
\end{equation}
is a minimal solution of the generalized backward stochastic
differential equation (GBSDE)
\begin{align}
\label{eq3.04} Y^{s,x}_t&=\varphi(X_T)
+\int^T_tf(\theta,X_{\theta},Y^{s,x}_{\theta},Z^{s,x}_{\theta})\,d\theta
+\int^T_tg(\theta,X_{\theta},Y^{s,x}_{\theta})\,dR_{s,\theta}\nonumber\\
&\quad-\int^T_tZ^{s,x}_{\theta}\,dB_{s,\theta},\quad
t\in[s,T],\quad P_{s,x}\mbox{-}a.s.,
\end{align}
where $B_{s,\cdot}$ is some standard Wiener process under
$P_{s,x}$ and $R$ is a continuous additive functional of
$\mathbb{X}$ corresponding to $\mu$ in the sense that
\begin{equation}
\label{eq1.06} E_{s,x}\int_{s}^{T}\eta(t,X_t)\,dR_t
=\int_s^T\!\!\int_{\mathbb{R}^d}\eta(t,y)p(s,x,t,y)\,d\mu(t,y)
\end{equation}
for any bounded measurable $\eta:Q_T\rightarrow[0,\infty)$.
Here $E_{s,x}$ denotes expectation with respect to $P_{s,x}$ and
$p$ is the transition density function of $(X,P_{s,x})$ (or,
equivalently, $p$ is the fundamental solution for $L_t$). From
(\ref{eq1.04}) it follows in particular that $u(s,x)=Y^{s,x}_s$,
so (\ref{eq3.04}) may be viewed as the Feynman-Kac formula for
solutions of (\ref{eqm.0}).

In \cite{EKPPQ} it is proved that viscosity solutions of the
Cauchy problem for semilinear parabolic equation in nondivergence
form with obstacle can be represented by solutions of some
reflected backward stochastic differential equations (RBSDEs). As
an application of results concerning (\ref{eqm.0}) we provide such
a representation in the case where the equation is in divergent
form and strong solutions in Sobolev spaces are considered. We
strengthen also known analytical results on homographic
approximation of solution of the obstacle problem.

Roughly speaking, the  obstacle problem consists in finding $u:Q_T
\rightarrow\mathbb{R}$ such that for given $\varphi, f$ as above
and $h:Q_T\rightarrow\mathbb{R}$,
\begin{equation}
\label{eq1.03} \left\{
\begin{array}{ll} \min(u-h,-\frac{\partial u}{\partial t}
-L_{t}u-f_u)=0 & \mbox{in }Q_T,\\
u(T)=\varphi & \mbox{on }\mathbb{R}^d,
\end{array}
\right.
\end{equation}
i.e. $u$ satisfies the prescribed terminal condition, takes values
above a given obstacle $h$,  satisfies inequality $\frac{\partial
u}{\partial t}+L_tu\leq -f_u$ in $Q_T$ and equation
$\frac{\partial u}{\partial t}+L_tu=-f_u$ on the set $\{u>h\}$.

In the case where $L_t$ is a non-divergent operator of the form
\[
L_{t}=\frac{1}{2}\sum_{i,j=1}^da^{ij} \frac{\partial^{2}}{\partial
x_{i}\partial x_{j}}+\sum_{i=1}^db^{i}\frac{\partial}{\partial
x_{i}}\,,
\]
problem (\ref{eq1.03}) has been investigated carefully in
\cite{EKPPQ} by using probabilistic methods. Let $X^{s,x}$ be a
solution of the It\^o equation
\[
dX^{s,x}_{t}=\sigma(t,X^{s,x}_t)\,dW_t+b(t,X^{s,x}_t)\,dt, \quad
X^{s,x}_s=x\quad (\sigma\sigma^*=a)
\]
associated with $L_t$. In \cite{EKPPQ} it is proved, that under
suitable assumptions on $a,b$ and the data $\varphi,f,h$, for each
$(s,x)\in Q_T$  there exists a unique solution
$(Y^{s,x},Z^{s,x},K^{s,x})$ of RBSDE with forward driving process
$X^{s,x}$, terminal condition $\varphi(X^{s,x}_T)$, coefficient
$f$ and obstacle $h(\cdot,X^{s,x}_{\cdot})$, and moreover, $u$
defined by the formula $u(s,x)=Y^{s,x}_s$, $(s,x)\in Q_T$ is a
unique viscosity solution of (\ref{eq1.03}) in the class of
functions satisfying the polynomial growth condition.

In the present paper we are interested in stochastic
representation of solutions of the obstacle problem with
divergence form operator in the framework of Sobolev spaces (for
the case of non-divergence form operator see \cite{Bally}, \cite{Ouknine}). The
advantage of using such a framework lies in the fact that it
allows to provide stochastic representation not only for $Y^{s,x}$
but also for $Z^{s,x}$ and $K^{s,x}$.

By the strong solution of the obstacle problem we understand a
pair $(u,\mu)$ consisting of a measurable function
$u:Q_T\rightarrow\mathbb{R}$ having some regularity properties and
a Radon measure $\mu$ on $Q_T$ such that
\begin{equation}
\label{eq1.4} \frac{\partial u}{\partial t}+L_tu=-f_u-\mu,\quad
u(T)=\varphi,\quad u\geq h,\quad \int_{Q_T}(u-h)\,d\mu=0
\end{equation}
(see Section \ref{sec3} for details).

Let $S_t=h(t,X_t)$, $t\in[s,T]$, and let
$(Y^{s,x},Z^{s,x},K^{s,x})$ be a solution of RBSDE
\begin{equation}
\label{eq1.6} \left\{
\begin{array}{l}
Y^{s,x}_t=\varphi(X_T)
+\int_t^Tf(\theta,X_\theta,Y^{s,x}_\theta,Z^{s,x}_\theta)\,d\theta
+K^{s,x}_T-K^{s,x}_t\medskip \\
\qquad\qquad\qquad\qquad\qquad\qquad\qquad\qquad
-\int_t^TZ^{s,x}_\theta\,dB_{s,\theta},\quad
t\in[s,T],\, P_{s,x}\mbox{-a.s.}\medskip\\
Y^{s,x}_t\ge S_t,\quad t\in[s,T],\medskip \\
K^{s,x}\mbox{  increasing, continuous, }K^{s,x}_s=0,\,\,
\int^T_s(Y^{s,x}_t-S_t)\,dK^{s,x}_t=0.
\end{array}
\right.
\end{equation}
We show that under mild conditions on $\varphi,f$ and $h$ there
exists a unique solution $(u,\mu)$ of (\ref{eq1.4}). Moreover, for
a.e. $(s,x)\in[0,T)\times{\mathbb{R}}^d$,
\begin{equation}
\label{eq1.5} u(t,X_t)=Y^{s,x}_t,\,t\in[s,T],\quad P_{s,x}
\mbox{-}a.s.,\qquad Z^{s,x}_t=\sigma\nabla u(t,X_t),\quad
\lambda\otimes P_{s,x}\ \mbox{-}a.s.
\end{equation}
and $K^{s,x}$ corresponds to $\mu$, i.e. (\ref{eq1.06}) with $R$
replaced by $K^{s,x}$ holds true. The correspondence between
$K^{s,x}$ and $\mu$ allows us to derive properties of $K^{s,x}$
from those of $\mu$ and vice versa.

Our proof of (\ref{eq1.5}) and the correspondence between
$K^{s,x}$ and $\mu$ is based on a general approximation result for
solutions of RBSDEs. The approximation we consider may be viewed
as an analogue of the well known in PDEs theory homographic
approximation for strong solutions of an obstacle problem (see
\cite{Palmeri}). Therefore we call it a stochastic homographic
approximation. Up to our knowledge, it is used here for the first
time in the context of RBSDEs.

By using the stochastic homographic approximation we prove also
that under mild regularity conditions on $h$ the measure $\mu$ is
absolutely continuous with respect to the Lebesgue measure
$\lambda$ on ${\mathbb{R}^d}$, and we get some information on the
density $d\mu/d\lambda$. This  provides information  on the
density of the control process $K^{s,x}$. It is worth pointing out
that the approximation provides additional information on the
control process also for general non-Markovian RBSDE with obstacle
being a general continuous semimartingale. For instance, it allows
to prove a stochastic version of the Lewy-Stampacchia inequality.

Our results on convergence of stochastic homographic
approximations to solutions of (\ref{eq1.6}) when combined with
(\ref{eq1.5}) prove convergence of homographic approximations of
solutions of (\ref{eq1.4}). In particular, we show that if
$\frac{\partial h}{\partial t}+L_th$ is a signed Radon measure on
$Q_T$ then under some assumptions on $\varphi,f$, the strong
solution $u$ of (\ref{eq1.4}) is a limit, in the space
${\mathbb{L}}_{2}(0,T;H^1_\varrho)\cap C([0,T],
{\mathbb{L}}_{2}({\mathbb{R}^d}))$, of maximal solutions of the
problem
\[
\frac{\partial u_{n}}{\partial t}
+L_{t}u_{n}=-f_{u_{n}}-\mu_n,\quad u_{n}(T)=\varphi
\]
with
\[
\mu_n=\frac{1}{1+n|u_{n}-h|}\,\Phi^-,\quad \Phi=\frac{\partial
h}{\partial t} +L_{t}h+f_{h}.
\]
This  strengthens considerably analytical results  which asserts
that $u$ is approximated by homographic approximations in
${\mathbb{L}}_{2,\varrho}(Q_T)$, while its gradient weakly in
${\mathbb{L}}_{2,\varrho}(Q_T)$ (see \cite{Palmeri}). Let us point
out also that contrary to \cite{Palmeri} our approximation is
direct in the sense that it does not require smoothing  the
functional $\Phi^-$. Furthermore, we prove  that $\{\mu_n\}$
converges to $\mu$ weakly and in the dual space to the space
${\mathcal{W}}_\varrho$ (see notation below).

Finally, let us mention that in case $b=0$ from  our stochastic
Lewy-Stampacchia inequality we get easily the Lewy-Stampacchia
inequality for solutions of the obstacle problem (\ref{eq1.4}).
\medskip

In the paper we adopt the following notation.

$Q_T=[0,T]\times{\mathbb{R}}^d$,
$Q_{\hat{T}}=[0,T)\times{\mathbb{R}}^d$,
$\check{Q}_{T}=(0,T)\times\mathbb{R}^{d}$,
$\nabla=(\frac{\partial}{\partial x_1},
\dots,\frac{\partial}{\partial x_d})$.

By $\mathcal{B}(D), \mathcal{B}_{b}(D), \mathcal{B}^{+}(D)$ we
denote  the set of Borel, bounded Borel, positive Borel functions
on $D$ respectively.
$C_{0}(D)$, $C_{0}^{\infty}(D)$, $C_{b}^{\infty}(D)$ are spaces of
all continuous functions with compact support in $D$, smooth
functions with compact support in $D$ and smooth functions on $D$
with bounded derivatives, respectively. We also write that $K\subset\subset D$
if $K$ is compact and included in $D$.

${\mathbb{L}}_{p}({\mathbb{R}}^d)$ (${\mathbb{L}}_{p}(Q_T)$) are
usual Banach spaces of measurable functions  on ${\mathbb{R}}^d$
(on $Q_T$) that are $p$-integrable. Let $\varrho$ be a positive
function on ${\mathbb{R}}^d$. By
${\mathbb{L}}_{p,\varrho}({\mathbb{R}}^d)$
(${\mathbb{L}}_{p,\varrho}(Q_T)$) we denote the space of functions
$u$ such that $u\varrho\in {\mathbb{L}}_p({\mathbb{R}}^d)$
($u\varrho\in{\mathbb{L}}_{p}(Q_T)$) equipped with the norm
$\|u\|_{p,\varrho}=\|u\varrho\|_p$
($\|u\|_{p,\varrho,T}=\|u\varrho\|_{p,T})$.
By $\langle\cdot,\cdot\rangle_{2,\varrho}$ we denote the inner product in
${\mathbb{L}}_{2,\varrho}({\mathbb{R}}^d)$  and by
$\langle\cdot,\cdot\rangle_{2,\varrho,T}$ the inner product in
${\mathbb{L}}_{2,\varrho}(Q_{T})$.

$H^1_{\varrho}$  is the Banach space consisting of all elements
$u$ of ${\mathbb{L}}_{2,\varrho}({\mathbb{R}}^d)$ having
generalized derivatives $\frac{\partial u}{\partial x_i}$,
$i=1,\dots,d$, in ${\mathbb{L}}_{2,\varrho}({\mathbb{R}}^d)$.
${\mathcal{W}}_\varrho$ is the subspace of ${\mathbb{L}}_2
(0,T;H^1_{\varrho})$ consisting of all elements $u$ such that
$\frac{\partial u}{\partial t}
\in{\mathbb{L}}_2(0,T;H^{-1}_{\varrho})$, where $H^{-1}_{\varrho}$
is the dual space to $H^{1}_{\varrho}$ (see \cite{Lions} for
details). By $\langle\cdot,\cdot\rangle_{\varrho,T}$ we denote
duality between ${\mathbb{L}}_2(0,T;H^{-1}_{\varrho})$ and
${\mathbb{L}}_2(0,T;H^{1}_{\varrho})$. $\mathcal{M}(Q_{T})$
$(\mathcal{M}^{+}(Q_{T}))$ denotes the space of Radon measures
(positive Radon measures) on $Q_{T}$. By $m$ we denote the
Lebesgue measure on $Q_{T}$.

By $C$ we  denote a general constant which may vary from line to
line but depends only on fixed parameters.

\nsubsection{Generalized BSDEs}
\label{sec2}

Let $\{B_{t}, 0\leq t\leq T\}$ be a $d$-dimensional standard
Wiener process  defined on some probability space $(\Omega,
\mathcal{F}, P)$ and let $\{\mathcal{F}_{t}\}_{t\in[0,T]}$ denote
the usual augmentation of the natural filtration generated  by
$B$.

Let $\xi$ be an $\mathcal{F}_{T}$-measurable random variable and
let $f:\Omega \times [0,T]\times
\mathbb{R}\times\mathbb{R}^{d}\rightarrow \mathbb{R}$, $g:\Omega
\times [0,T]\times \mathbb{R}\rightarrow \mathbb{R}$. We will need
the following assumptions.
\begin{enumerate}
\item[(A1)]
$\xi \in \mathbb{L}_{2}(\Omega,\mathcal{F}_{T},P)$,
\item[(A2)]
For every $(y,z)\in \mathbb{R}\times\mathbb{R}^{d}$ the processes
$f(\cdot,y,z)$, $g(\cdot,y)$ are predictable,
\item[(A3)]
$R$ is  increasing continuous process such that
$E|R_{T}|^2<\infty$,
\item[(A4)]
There exist $K>0$ and a  predictable process $\gamma$ such that
$E\int_{0}^T|\gamma_{t}|^{2}\,dt<\infty$ and
\[
|f(t,y,z)|\le K(|\gamma_{t}|+|y|+|z|),\quad P\mbox{-}a.s.
\]
for all $t\in [0,T], y\in\mathbb{R}, z\in\mathbb{R}^{d}$
\item[(A5)]
There exists $M>0$ such that $|g(t,\cdot)|\le M$, $t\in [0,T]$,
$P$-a.s.,
\item[(A6)]
$(y,z) \rightarrow f(t,y,z)$ is $P$-a.s. continuous for every
$t\in [0,T]$,
\item[(A7)]
$y\rightarrow g(t,y)$ is $P$-a.s. continuous for every $t\in
[0,T]$,
\item[(A6')]There is $L>0$ such that
\[
|f(t,y_1,z_1)-f(t,y_2,z_2)| \leq L(|y_1-y_2|+|z_1-z_2|),\quad
P\mbox{-a.s.}
\]
for every $y_1,y_2\in\mathbb{R},\,z_1,z_2\in\mathbb{R}^{d},
t\in[0,T]$,
\item[(A7')]There is $L>0$ such that
\[
|g(t,y_1)-g(t,y_2)| \leq L|y_1-y_2|,\quad  P\mbox{-a.s.}
\]
for every $y_1,y_2\in\mathbb{R}, t\in [0,T]$.
\end{enumerate}

Following \cite{PardouxZhangGRBSDE,RenXiaGRBSDE} we say that a
pair $(Y,Z)$ of $\{\mathcal{F}_{t}\}$-progressively measurable
processes on $[0,T]$ taking values in
$\mathbb{R}\times\mathbb{R}^{d}$ is a solution of the generalized
backward stochastic differential equation (GBSDE)
\begin{equation}
\label{eq2.02}
Y_{t}=\xi + \int_{t}^{T}f(s,Y_{s},Z_s)\,ds
+\int_{t}^{T}g(t,Y_{t})\,dR_{t} - \int_{t}^{T}(Z_{s},dB_{s}),\quad
t\in[0,T]
\end{equation}
if $E\sup_{0\leq t\leq T}|Y_{t}|^{2}<\infty$,
$E\int_{0}^{T}|Z_{t}|^{2}\,dt<\infty$ and (\ref{eq2.02}) is
satisfied $P$-a.s. If $(Y,Z)$ is a solution of (\ref{eq2.02}) such
that $Y_t\le\tilde Y_t$, $t\in[0,T]$, $P$-a.s. for any solution
$(\tilde Y,\tilde Z)$ of (\ref{eq2.02}), then it is called a
minimal solution of (\ref{eq2.02}).

Of course, if $R=0$ or $g=0$, then GBSDE
reduces to the usual  backward SDE with terminal condition $\xi$
and coefficient $f$.

The main purpose of the present section is to prove comparison
results for solutions of (\ref{eq2.02}), and, in consequence, to
prove that under (A1)--(A7) there is a  minimal solution to
(\ref{eq2.02}).

We begin with a priori estimates for solutions of (\ref{eq2.02})
and a "backward version" of the generalized Gronwall's lemma,
which in turn will be used to prove some comparison principle for
solutions of GBSDEs. Let us mention here that a priori estimates
and comparison results for solutions of GBSDEs are proved in
\cite{PardouxZhangGRBSDE} but under assumptions on $g$ not
suitable for our purposes (in \cite{PardouxZhangGRBSDE}
monotonicity of $g$ is required).

\begin{stw}\label{stwg.1}
Assume (A1)-(A5) and let $(Y,Z)$ be a solution of (\ref{eq2.02}).
Then there exists $C>0$ depending on $K,M,T$ such that
\[
E\sup_{0\le t\le T}|Y_{t}|^{2}+E\int_{0}^{T}|Z_{t}|^{2}\,dt \le
C(E|\xi|^2+E|R_{T}|^2+E\int_{0}^{T}|{\gamma}_{t}|^{2}\,dt).
\]
\end{stw}
\begin{dow}
By It\^o's formula, for every $t\in[0,T]$,
\[
|Y_{t}|^{2}+\int_{t}^{T}|Z_{\theta}|^{2}\,d\theta
=|\xi|^{2}+\int_{t}^{T}f(\theta,Y_{\theta},Z_{\theta})Y_{\theta}\,d\theta
+\int_{t}^{T}g(Y_{\theta})Y_{\theta}\,dR_{\theta}
-\int_{t}^{T}Z_{\theta}Y_{\theta}\,dB_{\theta}.
\]
Hence, by (A4) and (A5),
\begin{align}\label{eqg.2}
&|Y_{t}|^{2}+\frac12\int_{t}^{T}|Z_{\theta}|^{2}\,d\theta
\le|\xi|^{2}+2(K+K^2)
\int_{t}^{T}|Y_{\theta}|^2\,d\theta\nonumber\\
&\qquad+K\int_{t}^{T}|\gamma_{\theta}|^2\,d\theta
+M\int_{t}^{T}|Y_{\theta}|\,dR_{\theta}
-\int_{t}^{T}Z_{\theta}Y_{\theta}\,dB_{\theta}.
\end{align}
Taking expectation and using Gronwall's lemma yields
\begin{equation}
\label{eqg.3} E|Y_{t}|^{2}+E\int_{0}^{T}|Z_{\theta}|^{2}\,d\theta
\le C(E|\xi|^2+E\int_{0}^{T}|\gamma_{t}|^2\,dt
+E\int_{0}^{T}|Y_{\theta}|\,dR_{\theta}).
\end{equation}
Therefore taking supremum in (\ref{eqg.2}), using the
Burkholder-Davis-Gundy inequality  and (\ref{eqg.3}) we get
\begin{align*}
&E\sup_{0\le t\le T}|Y_{t}|^{2}
+E\int_{0}^{T}|Z_{\theta}|^{2}\,d\theta\\
&\qquad \le C(E\xi^{2}+E\int_{0}^{T}|\gamma_{t}|^{2}\,dt
+E\int_{0}^{T}|Y_{\theta}|dR_{\theta})\\
&\qquad\le
C(E|\xi|^{2}+E\int_{0}^{T}|\gamma_{t}|^{2}\,dt+E|R_{T}|^2)
+\frac12 E\sup_{0\le t\le T}|Y_{t}|^{2},
\end{align*}
which proves the proposition.
\end{dow}

\begin{lm}\label{lmg.1}
Let $Y$  be a continuous  decreasing process such that
$Y\ge0$ a.s. and  $EY_{0}<\infty$, and let $A$ be  an adapted
continuous increasing processes such that $0\le A_{T}\le a$ a.s.
for some $a>0$. If
\[
EY_{\tau}\le E\int_{\tau}^{T}Y_{s}\,dA_{s}
\]
for every stopping time $0\le \tau\le T$ then $Y=0$.
\end{lm}
\begin{dow}
Without lost of generality we may and will assume that $A$ is
strictly increasing. Put $\tau_{t}=\inf\{s\in [0,T]; A_{s}\ge t
\}\wedge T$. By the  change of variable formula  and assumptions,
\begin{align*}
EY_{\tau_{t}}&\le E\int^{T}_{\tau_{t}} Y_{s}\,dA_{s}
=E\int_{0}^{T}Y_{s}\mathbf{ 1}_{\{\tau_{t}\le s\le T\}}\,dA_{s}\\
&\le E\int_{0}^{\infty}
Y_{{\tau}_{u}}\mathbf{1}_{\{\tau_{t}\le\tau_{u}\le T\}}\,du
=\int_{t}^{A_{T}}EY_{\tau_{u}}\,du\le
\int_{t}^{a}EY_{\tau_{u}}\,du .
\end{align*}
and the result follows by classical Gronwall's lemma.
\end{dow}

\begin{tw}\label{thg.1}
Suppose that $\xi_{i},f_{i},g_{i},R^{i}$, $i=1,2$ satisfy
(A1)-(A5) and, in addition, $f_{1},g_{1}$ satisfy (A6'), (A7').
Let $(Y^{i},Z^{i})$ be a solution of (\ref{eq2.02}) with data
$\xi_{i}, f_{i}, g_{i},R^{i}$, $i=1,2$. If
\begin{enumerate}
\item [\rm(i)] $\xi_{1}\le \xi_{2}$, $P$-a.s.,
\item [\rm(ii)] $f_{1}(\cdot,y,z)\le f_{2}(\cdot,y,z)$, $dt\otimes
dP$-a.e. for every $(y,z)\in \mathbb{R}\times\mathbb{R}^{d}$,
\item [\rm(iii)]$g_{1}(t,Y^{2}_{t})\,dR^{1}_{t}
\le g_{2}(t,Y^{2}_{t})\,dR^{2}_{t}$, $P$-a.s.,
\item [\rm(iv)] $R^{1}_{T}\le a$, $P$-a.s. for some $a>0$,
\end{enumerate}
then $Y^{1}_{t}\le Y^{2}_{t}$, $t\in [0,T]$, $P$-a.s..
\end{tw}
\begin{dow}
Write $\xi=\xi_1-\xi_2$, $Y=Y^1-Y^2$ , $Z=Z^1-Z^2$. By It\^o's
formula and assumptions,
\begin{align*}
&|Y_{t}^{+}|^{2}+\int_{t}^{T}|Z_{\theta}|^{2}\mathbf{1}_{\{Y_{\theta}>0\}}\,
d\theta\\
&\quad =|\xi^{+}|^{2}+2\int_{t}^{T}
(f_{1}(\theta,Y^{1}_{\theta},Z^{1}_{\theta})
-f_{2}(\theta,Y^{2}_{\theta},Z^{2}_{\theta}))Y_{\theta}^{+}\,d\theta\\
&\qquad+\int_{t}^{T}g_{1}(\theta,Y^{1}_{\theta})
Y_{\theta}^{+}\,dR^{1}_{\theta}
-\int_{t}^{T}g_{2}(\theta,Y^{2}_{\theta})Y_{\theta}^{+}\,dR^{2}_{\theta}
-\int_{t}^{T}Z_{\theta}Y_{\theta}^{+}\,dB_{\theta}\\
&\quad\le C\int_{t}^{T}|Y_{\theta}^{+}|^{2}
+\alpha\int_{t}^{T}|Z_{\theta}|^{2}\mathbf{1}_{\{Y_{\theta}>0\}}\,d\theta
+\int_{t}^{T}(g_{1}(\theta,Y^{1}_{\theta})-g_{1}(\theta,Y^{2}_{\theta}))
Y_{\theta}^{+}\,dR^{1}_{\theta}\\
&\qquad -\int_{t}^{T}Z_{\theta}Y_{\theta}^{+}\,dB_{\theta}.
\end{align*}
Now, arguing as in the proof of Proposition \ref{stwg.1} we get
\begin{eqnarray*}
E\sup_{t\le s\le T}|Y_{s}^{+}|^{2}\le C\left(E
\int_{t}^{T}|Y_{\theta}^{+}|^{2}\,d\theta
+E\int_{t}^{T}|g_{1}(\theta,Y^{1}_{\theta})-g_{1}(\theta,Y^{2}_{\theta})|
Y_{\theta}^{+}\,dR^{1}_{\theta}\right).
\end{eqnarray*}
Hence,
\begin{eqnarray*}
E\sup_{t\le s\le T}|Y_{s}^{+}|^{2}\le CE\int_{t}^{T}\sup_{s\le
\theta\le T}|Y_{\theta}^{+}|^{2}(ds+dR^{1}_{s}).
\end{eqnarray*}
The same inequality we can get for every stopping time $\tau\le T$ instead of
$t$. By Lemma \ref{lmg.1} we get the result.
\end{dow}
\medskip

Let $\mathbb{Q}$ denote the set of rational numbers. The following
useful approximation result is proved in \cite{Lepeltier}.
\begin{lm}\label{lmm.1}
Let $f:\mathbb{R}^{d}\rightarrow \mathbb{R}$  be a continuous
function such that  $|f(x)|\le C(1+|x|)$,  $x\in \mathbb{R}^{d}$,
for some $C>0$. Set
$f_{n}(x)=\inf_{y\in\mathbb{Q}^{d}}\{f(y)+n|x-y|\}$,
$x\in{\mathbb{R}}^d$, $n\in{\mathbb{N}}$. Then
\begin{enumerate}
\item [\rm(a)] $|f_{n}(x)|\le C(1+|x|),\, x\in \mathbb{R}^{d},$
\item [\rm(b)] $f_{n}(x)\uparrow f(x),\, x\in \mathbb{R}^{d},$
\item [\rm(c)] if $x_{n}\rightarrow x$ then $f_{n}(x_{n})\rightarrow f(x),$
\item [\rm(d)] $f_{n}$ is Lipschitz continuous.
\end{enumerate}
\end{lm}

We will need also the following lemma.
\begin{lm}\label{lmm.2}
Let $f:\mathbb{R}^{d}\rightarrow[0,\infty)$  be a continuous
function. Then there is a sequence $\{f_{n}\}\subset
C_{0}^{\infty}(\mathbb{R}^{d})$ such that
\begin{enumerate}
\item [\rm(a)] $f_{n}(x)\uparrow f(x),\, x\in \mathbb{R}^{d},$
\item [\rm(b)] if $x_{n}\rightarrow x$ then $f_{n}(x_{n})\rightarrow f(x)$.
\end{enumerate}
\end{lm}
\begin{dow} By the Stone-Weierstrass theorem, for
every $\varepsilon>0$ there is $\tilde{f}_{\varepsilon}\in
C^{\infty}(Q_{T})$ such that
$\|\tilde{f}_{\varepsilon}-f\|_{\infty}\le \varepsilon$. Let
$\bar{f}_{n}=\tilde{f}_{4^{-n}}-2\cdot4^{-n}$. Then
$\bar{f}_{n}\le \bar{f}_{n+1}$, $n\ge1$, because $\bar{f}_{n}\le
f$ and $4^{-n}\le f-\bar{f}_{n}\le3\cdot4^{-n}$ for $n \ge1$.
Therefore the sequence $\{f_{n}=\eta_{n}\bar{f}_{n}\}$, where
$\{\eta_{n}\}\subset C_{0}^{\infty}(\mathbb{R}^d)$ is a sequence
of positive functions such that $\eta_{n}\uparrow1$ uniformly in
compacts subsets of $\mathbb{R}^d$ has the desired properties.
\end{dow}

\begin{stw}\label{prog.3}
If assumptions (A1)--(A7) are satisfied then there exists a
minimal solution of GBSDE (\ref{eq2.02}). Moreover, if
$\xi_{i},f_{i},g_{i},R^{i}$, $i=1,2$, satisfy assumptions
(i)--(iii) of Theorem \ref{thg.1} and the pairs $(Y^{i},Z^{i})$,
$i=1,2$, are minimal solutions of (\ref{eq2.02}) with data
$\xi_{i},f_{i},g_{i},R^{i}$, respectively, then $Y^{1}_{t}\le
Y^{2}_{t}$, $t\in[0,T]$, $P$-a.s..
\end{stw}
\begin{dow}
Let $f_n$ be the approximation of $f$ considered in Lemma
\ref{lmm.1} and let $g_n$ be the approximation of $g$ considered
in Lemma \ref{lmm.2}. From \cite{PardouxZhangGRBSDE} we know that
for each $n\in\mathbb{N}$  there exists a unique solution
$(Y^{n},Z^{n})$ of GBSDE
\begin{equation}
\label{eq2.04} Y^{n}_{t}=\xi +
\int_{t}^{T}f_{n}(s,Y^{n}_{s},Z^{n}_s)\,ds
+\int_{t}^{T}g_{n}(t,Y^{n}_{t})\,dR^{n}_{t} -
\int_{t}^{T}(Z^{n}_{s},dB_{s}),\quad t\in[0,T]
\end{equation}
with $R^n=R\wedge n$. By Theorem \ref{thg.1},  $\{Y^{n}\}$ is
increasing, and by Proposition \ref{stwg.1},
\begin{equation}\label{eqg.4}
E\sup_{0\le t\le T}|Y^{n}|^{2}_{t}
+E\int_{0}^{T}|Z^{n}_t|^{2}\,dt\le C\left(E|\xi|^{2}+E|R_{T}|^2
+E\int_{0}^{T}|{\gamma}_{t}|^2\,dt\right)
\end{equation}
for some $C$ not depending on $n$.  Therefore,
\[
E\int_{0}^{T}|Y^{n}_t-Y^{m}_t|^{2}\,dt
+E\int_{0}^{T}|Y^{n}_t-Y^{m}_t|\,dR_t\rightarrow0.
\]
By It\^o's formula,
\begin{align}\label{eqg.5}
&|Y^{n}_{t}-Y^{m}_{t}|^{2}
+\int_{t}^{T}|Z^{n}_{\theta}-Z^{m}_{\theta}|^{2}\,d\theta\nonumber\\
&\quad=\int_{t}^{T}(Y^{n}_{\theta}-Y^{m}_{\theta})
(f_{n}(\theta,Y^{n}_{\theta},Z^{n}_{\theta})
-f_{n}(\theta,Y^{n}_{\theta},Z^{n}_{\theta}))\,d\theta\nonumber\\
&\qquad+\int_{t}^{T}(Y^{n}_{\theta}-Y^{m}_{\theta})
g_{n}(\theta,Y^{n}_{\theta})\,dR^{n}_{\theta}
-\int_{t}^{T}(Y^{n}_{\theta}-Y^{m}_{\theta})g_{m}(\theta,Y^{m}_{\theta})\,
dR^{m}_{\theta}\nonumber\\
&\qquad+\int_{t}^{T} (Y^{n}_{\theta}-Y^{m}_{\theta})
(Z^{n}_{\theta}-Z^{m}_{\theta})\,dB_{\theta}.
\end{align}
From the above and (\ref{eqg.4}) we conclude that
\begin{align}\label{eqg.6}
&E|Y^{n}_{t}-Y^{m}_{t}|^{2}
+E\int_{0}^{T}|Z^{n}_t-Z^{m}_t|^{2}\,dt\nonumber\\
&\quad\le C\left(\left(E\int_{0}^{T}|Y^{n}_t-Y^{m}_t|^{2}\,
d\theta\right)^{1/2}
+E\int_{0}^{T}|Y^{n}_t-Y^{m}_t|\,dR_t\right)\equiv I^{n,m}.
\end{align}
Now taking supremum in (\ref{eqg.5}), using BDG inequality and
estimate (\ref{eqg.6}) we get
\[
E\sup_{0\le t\le T}|Y^{n}_{t}-Y^{m}_{t}|^{2}
+E\int_{0}^{T}|Z^{n}_t-Z^{m}_t|^{2}\,dt\le I^{n,m}.
\]
Since we know that $I^{n,m}\rightarrow0$ as
$n,m\rightarrow\infty$, passing to the limit in (\ref{eq2.04})
proves existence of a solution $(Y,Z)$ of (\ref{eq2.02}).
Furthermore, by Theorem \ref{thg.1}, if  $(\tilde{Y},\tilde{Z})$
is a solution of  (\ref{eq2.02}) then $Y^{n}_{t}\le
\tilde{Y}_{t}$, $t\in [0,T]$, $P$-a.s. for each $n\in\mathbb{N}$,
which implies that  $Y\le \tilde Y$.

To prove the second part of the theorem, we approximate
$(Y^{1},Z^{1})$ in the same manner as above. Let
$\{(Y^{1,n},Z^{1,n})\}$ denote the approximating sequence. By
Theorem \ref{thg.1},  $Y^{1,n}_{t}\le Y^{2}_{t}$, $t\in [0,T]$,
$P$-a.s. for $n\in\mathbb{N}$, which yields $Y^1\le Y^2$.
\end{dow}

\begin{uw}
\label{rem2.7} {\rm Under assumptions (A1)--(A7) there exists a
maximal solution of GBSDE (\ref{eq2.02}). This follows from the
fact that if $\bar f(t,y,z)=-f(t,-y,-z)$, $\bar g(t,y)=-g(t,-y)$,
and if $(\bar Y,\bar Z)$ is a solution of (\ref{eq2.02}) with
$\xi,f,g$ replaced by $-\xi,\bar f,\bar g$, then the pair $(-\bar
Y,-\bar Z)$ is a solution of (\ref{eq2.02}). Therefore, if $(Y,Z)$
is a minimal solution of (\ref{eq2.02}) with data $-\xi,\bar
f,\bar g$, then $(-Y,-Z)$ is a maximal solution of (\ref{eq2.02}).
}
\end{uw}

\nsubsection{Stochastic homographic approximation} \label{sec4}

In what follows $S$ denote a continuous
$\{\mathcal{F}_{t}\}$-progressively measurable real-valued process
$S$ on $[0,T]$ such that
\begin{enumerate}
\item[(A8)]
$S_{T}\leq \xi$ $P$-a.s. and $E\sup_{0\leq t\leq T}
|S^{+}_{t}|^{2}<\infty$.
\end{enumerate}

Recall that a triple $(Y,Z,K)$ of
$\{\mathcal{F}_{t}\}$-progressively measurable processes on
$[0,T]$ taking values in
$\mathbb{R}\times\mathbb{R}^{d}\times\mathbb{R}^{+}$ is a solution
of the  reflected backward stochastic differential equation
(RBSDE)
\begin{equation}
\label{eq2.03} \left\{
\begin{array}{l}
Y_{t}=\xi+\int_{t}^{T}f(s,Y_{s},Z_{s})ds + K_{T} - K_{t} -
\int_{t}^{T}(Z_{s},dB_{s}),\quad t\in[0,T],\medskip\\
Y_{t}\ge S_{t},\quad t\in[0,T],
\medskip \\
K\mbox{ is increasing, continuous, } K_0=0,
\int_{0}^{T}(Y_{t}-S_{t})\,dK_{t}=0
\end{array}
\right.
\end{equation}
if $E\sup_{0\leq t\leq T}|Y_{t}|^{2}<\infty$,
$E\int_{0}^{T}|Z_{t}|^{2}\,dt<\infty$, $E|K_T|^2<\infty$ and
(\ref{eq2.03}) is satisfied $P$-a.s.
\medskip

In \cite{EKPPQ} it is proved that if (A1), (A2), (A4), (A7') and
(A8) are satisfied then (\ref{eq2.03}) has a unique solution.

In the following  theorem we assume that  $S$ is a continuous
semimartingale  admitting the decomposition
\begin{equation}
\label{eq2.4} S_{t}=S_{T}+\int_{t}^{T}f(s,S_{s},\tilde{Z}_{s})\,ds
-(C_{T}-C_{t})+(R_{T}-R_{t})-\int_{t}^{T}\tilde{Z}_{s}\,dB_{s},
\end{equation}
where $\tilde{Z}$ is an $\{\mathcal{F}_{t}\}$-adapted
square-integrable process and $C,R$ are continuous
$\{\mathcal{F}_{t}\}$-adapted square-integrable increasing
processes.

\begin{uw}
{\rm If  $S$ is a semimartingale with the decomposition
\[
S_t=S_0+M_t+U_t,\quad t\in[0,T],
\]
where $M$ is square-integrable martingale on $[0,T]$ and $U$ is an
adapted process of square-integrable  variation on $[0,T]$, then
it admits decomposition of the form (\ref{eq2.4}). Indeed, by the
representation theorem for martingales,  there is a progressively
measurable process $\tilde{Z}$ such that
$E\int^T_0|\tilde{Z}_s|^2\,ds<\infty$ and
\[
M_t=\int^t_0\tilde{Z}_s\,dB_s,\quad t\in[0,T].
\]
Moreover, since $U$ is a finite variation process, there exist
increasing processes $U^+$, $U^-$ such that
$U_{t}=U_{t}^{+}-U_{t}^{-}$, $t\in[0,T]$. Therefore putting
\[
C_t =\int^t_0(f(s,S_s,\tilde{Z}_s))^+\,ds+U^+_t, \quad
R_t=\int^t_0(f(s,S_s,\tilde{Z}_s))^-\,ds+U^-_t,\quad t\in[0,T]
\]
yields (\ref{eq2.4}).}
\end{uw}

\begin{tw}
\label{tw2.1} Assume (A1),(A2),(A4),(A7') and (A8). Let $S$ be of
the form (\ref{eq2.4}) and for $n\in\mathbb{N}$ let $(Y^n,Z^n)$ be
a maximal solution of the following GBSDE
\begin{equation}
\label{eq2.1}
Y^{n}_{t}=\xi+\int_{t}^{T}f(s,Y^{n}_{s},Z^{n}_{s})\,ds
+K^n_T-K^n_t-\int_{t}^{T}Z^{n}_{s}\,dB_{s},\quad t\in[0,T],
\end{equation}
where
\begin{equation}
\label{eq2.2}K^{n}_{t}=\int_{0}^{t}\alpha^n_s\,dR_s,\quad
\alpha^n_t=\frac{1}{1+n|Y^{n}_{t}-S_{t}|}.
\end{equation}
Then $Y^{n}_{t}\ge S_{t}$, $t\in[s,T]$ a.s. for each
$n\in\mathbb{N}$, $Y^{n}_{t}\downarrow Y_{t}$ a.s. for every
$t\in[0,T]$ and
\begin{equation}
\label{eq2.3} E\sup_{t\in[0,T]}|Y^{n}_{t}-Y_{t}|^{2}
+E\int_{0}^{T}|Z^{n}_{t}-Z_{t}|^{2}\,dt
+E\sup_{t\in[0,T]}|K^{n}_{t}-K_{t}|^{2}\rightarrow0,
\end{equation}
where $(Y,Z,K)$ is a solution of  (\ref{eq2.03}).
\end{tw}
\begin{dow}
Let $(\bar{Y}_{t}^{n},\bar{Z}_{t}^{n})$ be a solution of
(\ref{eq2.1}) with $\bar{\alpha}_{t}^{n}
=1-\frac{n(\bar{Y}_{t}^{n}-S_{t})}{1+n|\bar{Y}_{t}^{n}-S_{t}|}$ in
place of $\alpha^n$. By It\^{o}'s formula, for every $t\in[0,T]$
we have
\begin{align*}
|(\bar{Y}^{n}_{t}-S_{t})^{-}|^{2}&=|(\bar{Y}^{n}_{T}-S_{T})^{-}|^{2}
-2\int_{t}^{T}\frac{n|(\bar{Y}^{n}_{s}-S_{s})^{-}|^{2}}
{1+n|\bar{Y}^{n}_{s}-S_{s}|}\,dR_{s}\\
&\quad-2\int_{t}^{T}(\bar{Y}^{n}_{t}-S_{t})^{-}\,dC_{t}
-2\int_{t}^{T}(\bar{Y}^{n}_{t}-S_{t})^{-}
(\bar{Z}^{n}_{s}-\tilde{Z}_{s})\,dB_{s}\\
&\quad-\int_{t}^{T}\mathbf{1}_{\{\bar{Y}^{n}_{t}-S_{t}\leq0\}}
|\bar{Z}^{n}_{s}-\widetilde{Z}_{s}|^{2}\,ds\\
&\quad+2\int_{t}^{T}(f(s,S_{s},\tilde{Z}_{s})
-f(s,\bar{Y}^{n}_{s},\bar{Z}^{n}_{s}))(\bar{Y}^{n}_{s}-S_{s})^{-}\,ds.
\end{align*}
From this we obtain
\[
E|(\bar{Y}^{n}_{t}-S_{t})^{-}|^{2}\leq
CE\int_{t}^{T}\left|(\bar{Y}^{n}_{s}-S_{s})^{-}\right|^{2}\,ds,
\]
which implies that $S\le\bar{Y}^{n}$. From this we see that
$(\bar{Y}^{n},\bar{Z}^{n})$ is a solutions of  (\ref{eq2.1}). From
maximality of $Y^{n}$ we have that $S\le \bar{Y}^{n}\le Y^{n}$ and
we get (i). Using It\^o's formula, the Burkholder-Davis-Gundy
inequality and standard estimates we get
\[
E\sup_{0\le t\le T} |Y^{n}_t|^2+E\int^T_0|Z^n_t|^2\,dt \le
CE\left(|\xi|^2+\int_0^T|\gamma_{t}|^2\,dt
+\int_0^T\alpha^n_t|Y_t^n|\,dR_t \right).
\]
It follows from the form of equation (\ref{eq2.1}) and Proposition
\ref{prog.3} that $Y^{n}_{t}\geq Y^{n+1}_{t}$, $t\in[0,T]$,
$P$-a.s., $n\in\mathbb{N}$. Hence
\begin{align}
\label{eqq2.6}
&E\sup_{0\le t\le T}|Y^{n}_t|^2+E\int^T_0|Z^n_t|^2\,dt\nonumber\\
&\qquad\le
CE\left(|\xi|^2+\int_0^T|\gamma_{t}|^2\,dt+\int_0^T\alpha^n_t|Y_t^1|\,dR_t
+\int_0^T\alpha^n_t|S_t|\,dR_t\right).
\end{align}
Using once again It\^o's formula we get
\begin{eqnarray*}
&&E|Y^{n}_{t}-Y_{t}^{m}|^{2}+E\int_{t}^{T}|Z^{n}_{t}-Z_{t}^{m}|^{2}\,dt\\
&&\qquad= -2E\int_{t}^{T}(Y^{n}_{s}-Y_{s}^{m})
(f(s,Y^{m}_{s},Z^{m}_{s})-f(s,Y^{n}_{s},Z^{n}_{s}))\,ds\\
&&\qquad\quad-E\int_{t}^{T}(Y^{n}_{s}-Y_{s}^{m})
(\alpha^{m}_{s}-\alpha_{s}^{n})\,dR_{s}\\
&&\qquad \leq\frac{1}{2}E\int_{t}^{T}|Z^{n}_{s}-Z_{s}^{m}|^{2}\,ds
+CE\int_{t}^{T}|Y^{n}_{s}-Y_{s}^{m}|^{2}\,ds\\
&&\qquad\quad+E\int_{0}^{T}
|\alpha^{m}_{s}-\alpha_{s}^{n}||Y^{n}_{s}-Y_{s}^{m}|\,dR_{s}.
\end{eqnarray*}
By the above and Gronwall's lemma,
\begin{equation}
\label{eq2.8}
E|Y^{n}_{t}-Y^{m}_{t}|^{2}+\int_{0}^{T}|Z^{n}_{s}-Z^{m}_{s}|^{2}\,ds
\leq CE\int_{0}^{T}|\alpha^{m}_{s}-\alpha^{n}_{s}|
|Y^{n}_{s}-Y^{m}_{s}|\,dR_{s}.
\end{equation}
From the monotonicity of $\{Y^{n}\}$ there is a process $\bar{Y}$
such that $Y^{n}_{t}\searrow \bar{Y}_{t}$, $t\in[0,T]$. From this
and (\ref{eqq2.6}) we conclude that
$E\int_{0}^{T}|Y^{n}_{t}-Y^{m}_{t}|^{2}\,dR_{t}\rightarrow0$.
Hence, by (\ref{eq2.8}),
$E\int_{0}^{T}|Z^{n}_{t}-Z^{m}_{t}|^{2}dt\rightarrow0$ as
$n,m\rightarrow\infty$. Using the Burkholder-Davis-Gundy
inequality we conclude from the above that
\[
E\sup_{t\in[0,T]}|Y^{n}_{t}-Y^{m}_{t}|^{2}+E\int_{0}^{T}
|Z^{n}_{t}-Z^{m}_{t}|dt\leq
C\int_{0}^{T}|\alpha^{n}_{s}-\alpha^{m}_{s}|
|Y^{n}_{s}-Y^{m}_{s}|dR_{s}\rightarrow0
\]
as $n,m\rightarrow\infty$, and hence, by (\ref{eq2.1}), that
$E\sup_{t}|K^{n}_{t}-K^{m}_{t}|^{2}\rightarrow0$ as
$n,m\rightarrow\infty$. This implies that there is a triple
$(\bar{Y},\bar{Z},\bar{K})$ such that $\bar{Y}$ is continuous,
$\bar{K}$ continuous and increasing, satisfying
\[
E\sup_{t\in[0,T]}(Y^{n}_{t}-\bar{Y}_{t})^{2}
+E\int_{0}^{T}|Z^{n}_{t}-\bar{Z}_{t}|^{2}\,dt
+E\sup_{t\in[0,T]}|K^{n}_{t}-\bar{K}_{t}|^{2}\rightarrow0.
\]
From this we obtain that
\[
\int_{0}^{T}\mathbf{1}_{\{\bar{Y}_{s}-S_{s}>0\}}(Y^{n}_{s}-S_{s})\,dK^{n}_{s}
\rightarrow \int_{0}^{T}\mathbf{1}_{\{\bar{Y}_{s}-S_{s}>0\}}
(\bar{Y}_{s}-S_{s})\,d\bar{K}_{s}
=\int_{0}^{T}(\bar{Y}_{s}-S_{s})\,d\bar{K}_{s}.
\]
$P$-a.s.. On the other hand,
\[
\int_{0}^{T}\mathbf{1}_{\{\bar{Y}_{s}-S_{s}>0\}}
(Y^{n}_{s}-S_{s})\,d{K}^{n}_{s}
=\int_{0}^{T}\mathbf{1}_{\{\bar{Y}_{s}-S_{s}>0\}}
\frac{(Y^{n}_{s}-S_{s})}{1+n(Y^{n}_{s}-S_{s})}\,dR_{s}\rightarrow0.
\]
Accordingly, $\int_{0}^{T}(\bar{Y}_{s}-S_{s})\,d\bar{K}_{s}=0$,
$P\mbox{-a.s.}$. Therefore, by  uniqueness of solutions of RBSDEs,
$(\bar{Y},\bar{Z},\bar{K})=(Y,Z,K)$
\end{dow}
\medskip

The following corollary may be viewed as a stochastic version of
the Lewy-Stam\-pacchia inequality (see \cite{Donati,Palmeri} and
Remark \ref{uw5.8}).
\begin{wn}\label{col2.2}
Under assumptions of Theorem \ref{tw2.1},
\begin{equation}
\label{eq2.11} 0\le dK_{t}\le\mathbf{1}_{\{Y_{t}=S_{t}\}}\,dR_{t}.
\end{equation}
\end{wn}
\begin{dow}
Follows from (\ref{eq2.2}), (\ref{eq2.3}).
\end{dow}

\begin{uw}\label{rem2.3}
{\rm If $S$ is an It\^o process of the form
\[
S_{t}=S_{0}+\int_{0}^{t}\tilde{Z}_s\,dB_s+\int_{0}^{t}U_{s}\,ds,\quad
t\in[0,T],
\]
where $U,\tilde{Z}$ are progressively measurable processes such
that $E\int^T_0(|U_t|^{2}+|\tilde{Z}_t|^{2})\,dt<\infty$, then by
\cite[Remark 4.3]{EKPPQ},
\begin{equation}
\label{eq2.9} K_{t}=\int^{t}_{0}\mathbf{1}_{\{Y_{s}=S_{s}\}}
\alpha_{s}(f(s,S_{s},\tilde{Z}_{s})+U_{s})^{-}\,ds
\end{equation}
for some progressively measurable process $\alpha$ with values in
$[0,1]$. In fact, from every subsequence $\{n'\}$ we may choose a
further subsequence $\{n''\}$ such that
$\alpha^{n''}\rightarrow\alpha$ weakly in
$\mathbb{L}_2((0,T)\times\Omega;dt\otimes dP)$, which provides
some additional information on $\alpha$. To see that $\alpha$ can
be approximated by $\alpha_n$, let us first observe that $S$ may
be written in the form
\[
S_{t}=S_{T}+\int_{t}^{T} f(s,S_{s},Z_{s})\,ds-(\tilde C_T-\tilde
C_t)+(\tilde R_T-\tilde R_t)-\int_{t}^{T}Z_{s}\,dB_s, \quad
t\in[0,T],
\]
where
\[
\tilde
C_{t}=\int_{0}^{t}(f(s,S_{s},\tilde{Z}_{s})+U_{s})^{+}\,ds,\quad
\tilde R_{t}=\int_{0}^{t}(f(s,S_{s},\tilde{Z}_{s})+U_{s})^{-}\,ds.
\]
Therefore, by Theorem \ref{tw2.1}, if $Y^n,Z^n,\alpha^n$ are
defined by (\ref{eq2.1}), (\ref{eq2.2}) and $K^n$ is defined by
the formula
\[
K^{n}_{t}=\int_{0}^{t}
\alpha^{n}_{s}(f(s,S_{s},\tilde{Z}_{s})+U_{s})^{-}\,ds,
\]
then (\ref{eq2.3}) holds true. Since $\alpha^{n}$ are uniformly
bounded, there is a subsequence $\{n'\}$ such that
$\alpha^{n'}\rightarrow\bar\alpha$ weakly in
$\mathbb{L}_2((0,T)\times\Omega;dt\otimes dP)$. Since
$\alpha^{n}_{s}(f(s,S_{s},\tilde{Z}_{s})+U_{s})^{-}$ are uniformly
bounded in $\mathbb{L}_2((0,T)\times\Omega;dt\otimes dP)$ as well,
there is $\{n''\}\subset\{n'\}$ such that
\[
\alpha^{n''}_{t}(f(t,S_{t},\tilde{Z}_{t})+U_{t})^{-}\rightarrow
\mathbf{1}_{\{Y_{t}=S_{t}\}}
\alpha_{t}(f(t,S_{t},\tilde{Z}_{t})+U_{t})^{-},\quad
\alpha^{n''}\rightarrow \bar{\alpha}
\]
weakly in $\mathbb{L}_2((0,T)\times\Omega;dt\otimes dP)$. From
this we conclude that $\bar{\alpha}=\alpha$ on the set
$\{\mathbf{1}_{\{Y_{t}=S_{t}\}}
\alpha_{s}(f(s,S_{s},\tilde{Z}_{s})+U_{s})^{-}>0\}$. }
\end{uw}

\begin{uw}
{\rm Analysis of the proof of Theorem \ref{tw2.1} shows that the
assumption that $S$ is continuous is superfluous. What we really
need is continuity  of the process $R$. This is related to the
fact, that if $C$ is c\`adl\`ag and $R$ is continuous then $S$ has
only downward jumps going backward in time (see \cite{H}).}
\end{uw}

\nsubsection{Semilinear parabolic equations with measure data}
\label{sec3}

Let $\varphi:\mathbb{R}^d\rightarrow\mathbb{R}$,
$f:Q_T\times\mathbb{R}\times\mathbb{R}^{d} \rightarrow\mathbb{R}$,
$g:Q_T\times\mathbb{R}\rightarrow\mathbb{R}$ be measurable
functions. In this section we are concerned
with existence and stochastic representation of a solution of the
problem (\ref{eqm.0})

Let $\Omega=C([0,T],\mathbb{R}^d)$ denote the space of continuous
$\mathbb{R}^d$-valued functions on $[0,T]$ equipped with the
topology of uniform convergence and let $X$ be a canonical process
on $\Omega$. It is known that for an operator $L_t$ defined by
(\ref{eq1.3}) with $a$ and $b$ satisfying (\ref{eq1.1}) one can
construct a weak fundamental solution $p(s,x,t,y)$ for $L_t$ and
then a Markov family $\mathbb{X}=\{(X,P_{s,x});(s,x)\in[0,T)\times
\mathbb{R}^d\}$ for which $p$ is the transition density function,
i.e.
\[
P_{s,x}(X_t=x;0\leq t\leq s)=1,\quad
P_{s,x}(X_t\in\Gamma)=\int_{\Gamma}p(s,x,t,y)\,dy,\quad t\in(s,T]
\]
for any $\Gamma$ in a Borel $\sigma$-field $\mathcal{B}$ of
$\mathbb{R}^d$  (see \cite{Ro1,Stroock.DIFF}).

Set ${\mathcal{F}}^s_t=\sigma(X_u,u\in[s,t]),
\bar{\mathcal{F}}^{s}_{t}=\sigma(X_u,u\in[T+s-t,T])$ and define
${\mathcal{G}}$ as the completion of ${\mathcal{F}}^s_T$ with
respect to the family ${\mathcal{P}}=\{P_{s,\mu}:\mu$ is a
probability measure on ${\mathcal{B}}({\mathbb{R}}^d)$\}, where
$P_{s,\mu}(\cdot)=\int_{{\mathbb{R}}^d}P_{s,x}(\cdot)\,\mu(dx)$,
and define ${\mathcal{G}}^s_t$ ($\bar{\mathcal{G}}^s_t)$ as the
completion of ${\mathcal{F}}^s_t$ ($\bar{\mathcal{F}}^s_t)$ in
${\mathcal{G}}$ with respect to ${\mathcal{P}}$.

We will say that a family $A=\{A_{s,t}, 0\le s\le t\le T\}$ of
random variables is an additive functional (AF) of ${\mathbb{X}}$
if $A_{s,t}$ is ${\mathcal{G}}_{t}^{s}$-measurable for every $0\le
s\le t\le T$ and $P_{s,x}(A_{s,t}=A_{s,u}+A_{u,t}, s\le u\le t\le
T)=1$ for a.e. $(s,x)\in Q_{\hat{T}}$. If, in addition,
$A_{s,\cdot}$ has $P_{s,x}$-almost all continuous trajectories for
a.e. $(s,x)\in Q_{\hat{T}}$, then  $A$ is called a continuous AF
(CAF), and if  $A_{s,\cdot}$ is an increasing process under
$P_{s,x}$ for  a.e. $(s,x)\in Q_{\hat{T}}$, it  is called an
increasing AF. If $M$ is an AF such that for a.e. $(s,x)\in
Q_{\hat{T}}$, $E_{s,x}|M_{s,t}|^2<\infty$ and $E_{s,x}M_{s,t}=0$
for $t\in[s,T]$ it is called a martingale AF (MAF). Finally, we
say that $A$ is an AF (CAF, increasing AF, MAF) in the strict
sense if the corresponding property holds for every $(s,x)\in
Q_{\hat{T}}$.

Now we recall some known facts about  functionals in
${\mathbb{L}}_{2}(0,T;H^{-1}_{\varrho})$ (for details  see, e.g.,
 \cite{Evans}, \cite{Kl}). Here and in what follows we will assume that
$\varrho(x)=(1+|x|^{2})^{-\alpha}$, $x\in\mathbb{R}^d$, for some
$\alpha\ge0$ and $\int_{\mathbb{R}^{d}}\varrho(x)\,dx<\infty$.

It is  known that  if $\Phi\in\mathbb{L}_2(0,T;H^{-1}_{\varrho})$
then $\mu=f-\dyw{\bar{f}}$ for some $f,\bar{f}=(\bar
f_1,\dots,\bar f_d)\in \mathbb{L}_{2,\varrho}(Q_{T})$, i.e.
\begin{equation}\label{eqm.0.2}
\Phi(\eta)=\langle f,\eta \rangle_{2,\varrho,T} +\langle
\bar{f},\nabla\eta\rangle_{2,\varrho,T},\quad \eta\in
\mathbb{L}_{2}(0,T;H^1_{\varrho}).
\end{equation}
Let $\Phi\in{\mathbb{L}}_{2}(0,T;H^{-1}_{\varrho})$. We say that
$\Phi\in\mathcal{M}^{+}(Q_{T})$ if there is a measure $\mu\in
\mathcal{M}^{+}(Q_{T})$ such that
\begin{equation}
\label{RM}
\Phi(\eta)=\int_{Q_{T}}\eta\,d\mu
\end{equation}
for every $\eta\in C_{0}^{\infty}(Q_{T})$. Let us note that the
measure $\mu$ has the property that
$\mu({\{t\}\times{\mathbb{R}}^{d}})=0$ for every $t\in[0,T]$ (see
\cite{Kl}).

Let us consider the Dirichlet form $(\mathcal{E},\mathcal{F})$ on
$\mathbb{L}_{2,\varrho}(Q_{T})$ with
$\mathcal{F}=\mathbb{L}_{2}(0,T;H^{1}_{\varrho})$ defined by the
formula
\begin{equation*}
\mathcal{E}(u,v)=\langle\nabla u,\nabla
v\rangle_{2,\varrho,T}\,,\quad u,v\in\mathcal{F}.
\end{equation*}
It is easy to check that $\mathcal{E}$ is regular and
$C_{0}^{\infty}(Q_{T})$ is its core. With the form
$(\mathcal{E},\mathcal{F})$ we may associate a Choquet capacity
$\mbox{Cap}:2^{Q_{T}}\rightarrow[0,\infty]$ as follows. Let
$\mathcal{O}$ denote the family of all open subsets of $Q_{T}$.
For $A\in \mathcal{O}$ we put
\begin{equation*}
\mbox{Cap}(A)=\inf_{u\in\mathcal{L}_{A}}\mathcal{E}_{1}(u,u),
\end{equation*}
where $\mathcal{L}_{A}=\{u\in\mathcal{F}; u\ge 1\mbox{ a.e. on }
A\}$ and $\mathcal{E}_{1}(u,v)=\mathcal{E}(u,v)+\langle
u,v\rangle_{2,\varrho,T}.$ For $A\subset Q_{T}$ we put
\[
\mbox{Cap}(A)=\inf_{B\in\mathcal{O}, A\subset B} Cap(B).
\]
By \cite[Theorem 2.1.5]{Fukushima}, for every $A\subset Q_{T}$
there exists a unique
$e_{A}\in\bar{\mathcal{L}}_{A}=\{u\in\mathcal{F};u\ge 1\mbox{
Cap-q.e. on }A\}$ such that $Cap(B)=\mathcal{E}_{1}(e_{A},e_{A})$. Since every functional in
$\mathbb{L}_{2}(0,T;H^{-1}_{\varrho})$ is of the form
(\ref{eqm.0.2}), it follows that the measure $\mu$ is of finite
energy integral (see Section 2.2 in \cite{Fukushima}) and, by
\cite[Lemma 2.2.3]{Fukushima}, $\mu\ll\mbox{Cap}$. Moreover, since
every $\eta\in\mathcal{F}$ has Cap-quasi continuous version,
repeating arguments of  the proof of \cite[Theorem
2.2.2]{Fukushima} we can extend formula (\ref{RM}) to all
$\eta\in\mathcal{F}$.  In particular, given $\alpha\in
\mathcal{B}_{b}(Q_{T})$ and
$\Phi\in{\mathbb{L}}_2(0,T;H^{-1}_\varrho)\cap\mathcal{M}^{+}(Q_{T})$
we may define $\alpha\Phi\in{\mathbb{L}}_2(0,T;H^{-1}_\varrho)$ by
the formula
\[
\alpha\Phi(f)=\Phi(\alpha f)=\int_{Q_T}\alpha f\,d\mu,\quad f\in
{\mathbb{L}}_2(0,T;H^1_\varrho),
\]
where $\mu$ is the measure associated with $\Phi$ in the sense of
(\ref{RM}).

Let us now consider parabolic capacity naturally related to the
space $\mathcal{W}_{\varrho}$. We define the parabolic capacity of
the Borel set $B\subset Q_{\hat{T}}$ as follows
\[
\mbox{cap}(B)=P_{m}(\{\omega:(t,X_{t}(\omega))\in B\mbox{ for some
}t\in[0,T]\}),
\]
where $m$ is the Lebesgue measure on $Q_{T}$ and
$P_{m}(\Gamma)=\int_{Q_{\hat{T}}}P_{s,x}(\Gamma)\,ds\,dx$ for
$\Gamma\in\mathcal{G}$. We say that $u\in \mathcal{B}(Q_{T})$ is
cap-quasi continuous if $[s,T]\ni t\mapsto u(t,X_{t})$ is a
$P_{s,x}$-a.s. continuous process for a.e. $(s,x)\in Q_{T}$. It is
known (see \cite{Kl,Oshima}) that every $\eta\in
\mathcal{W}_{\varrho}$ has a cap-quasi continuous version. In what
follows we will always consider cap-quasi continuous versions of
elements of $\mathcal{W}_{\varrho}$.

From \cite[Theorem 2.1.4]{Fukushima} it follows that if
$u,\bar{u}$ are Cap-quasi continuous and $u=\bar{u}$ a.e. then
they are equal Cap-quasi everywhere. The same property hold for
parabolic capacity.
\begin{stw}
If $u,\bar{u}\in \mathbb{L}_{2,\varrho}(Q_{T})$ are cap-quasi
continuous and $u=\bar{u}$ a.e. then $u=\bar{u}$ cap-quasi
everywhere.
\end{stw}
\begin{dow}
Suppose that $\mbox{cap}(\{u\neq\bar{u}\})>0$. Then there exists
$A\subset Q_{\hat{T}}$ such that $m(A)>0$ and for every $(s,x)\in
A$,
\[
P_{s,x}(\{\omega:(t,X_{t})\in\{u\neq\bar{u}\}\mbox{ for some }
t\in[s,T]\})> 0
\]
Since the processes $t\mapsto u(t,X_t)$, $t\mapsto\bar u(t,X_t)$
have continuous trajectories,
\begin{align}
\label{eq5.3} 0<E_{s,x}\int_{s}^{T}|u-\bar{u}|^{2}(t,X_{t})\,dt
=\int_{s}^{T}\int_{\mathbb{R}^{d}}
|u-\bar{u}|^{2}p(s,x,\theta,y)\,d\theta\,dy.
\end{align}
Therefore, since $m(A)>0$,
\begin{align*}
0&<\int_{0}^{T}\int_{\mathbb{R}^{d}}
\left(E_{s,x}\int_{s}^{T}|u-\bar{u}|^{2}(t,X_{t})\,dt\right)
\varrho^{2}(x)\,ds\\
&\quad\le C\int_{0}^{T}\int_{\mathbb{R}^{d}}
|u-\bar{u}|^{2}(\theta,y)\varrho^{2}(y)\,d\theta\, dy,
\end{align*}
the last inequality being a consequence of \cite[Proposition
4.1]{Kl}. Since (\ref{eq5.3}) contradicts the assumption that
$u=\bar{u}$ a.e., the proposition is proved.
\end{dow}

\begin{uw}
{\rm  If $\eta\in \mathcal{W}_{\varrho}$ then from \cite[Appendix
A.2]{DPP} it follows that there exists $\{\eta_{n}\}\subset
C_{0}^{\infty}(Q_{T})$ such that $\eta_{n}\rightarrow \eta$ in
$\mathcal{W}_{\varrho}$. By \cite[Corollary 3.4]{Kl} there exists
a subsequence (still denoted by $\{n\}$) such that
$\eta_{n}\rightarrow\bar{\eta}$ cap-q.e., where $\bar{\eta}$ is
cap-quasi continuous version of $\eta$. On the other hand,
$\eta_{n}\rightarrow\eta$ in $\mathcal{E}_{1}$ so by \cite[Theorem
2.1.4]{Fukushima} there exists a subsequence (still denoted by
$\{n\}$) such that $\eta_{n}\rightarrow\tilde{\eta}$ Cap-q.e.,
where $\tilde{\eta}$ is Cap-quasi continuous version of $\eta$.
From this we conclude that $\int \bar{\eta}\,d\mu=\int
\tilde{\eta}\,d\mu$ for $\mu\in
\mathbb{L}_{2}(0,T;H^{-1}_{\varrho})\cap\mathcal{M}(Q_{T})$. }
\end{uw}

Let $\mu$ be a positive Radon measure on $Q_T$ and let  $K$ be an
increasing CAF of $\mathbb{X}$.  We will say that $\mu$
corresponds to $K$ or $K$ corresponds to  $\mu$  (and write
$\mu\sim K$) if
\begin{equation}
\label{eq3.01} E_{s,x}\int_{s}^{T}\eta(t,X_t)\,dK_{s,t}
=\int_{Q_{sT}}\eta(t,y)p(s,x,t,y)\,d\mu(t,y)
\end{equation}
for every $\eta\in\mathcal{B}^{+}(Q_{T})$ and a.e. $(s,x)\in Q_{T}$.

Observe that if $\mu$ corresponds to some increasing CAF of
$\mathbb{X}$, then $\mu\ll cap$ since $p>0$. Note also that
from \cite[Corollary 3.5]{Kl} it follows that every
$\mu\in\mathbb{L}_{2}(0,T;H^{-1}_{\varrho})\cap\mathcal{M}^{+}(Q_{T})$
has a corresponding CAF of $\mathbb{X}$.

Now we prove some properties of the Laplace transform of
time-inhomogeneous  additive functionals. Analogous properties for
time-homogenuous additive functionals are to be found for instance
in \cite[Chapter X]{Revuz}.

Let $A$ be an additive functional of $\mathbb{X}$ and let
$\alpha\ge0$. The function
\[
U^{\alpha}_{A}(s,x) =E_{s,x}\int_{s}^{T}e^{-\alpha
(t-s)}\,dA_{s,t},\quad (s,x)\in Q_T
\]
is called the Laplace transform of the AF $A$ or the
$\alpha$-potential of $A$. If $f\in \mathcal{B}_{b}(Q_{T})$ and
$f\cdot A$ is the functional defined by $(f\cdot A)_{s,t}
=\int^t_sf(\theta,X_{\theta})\,d\theta$, $0\le s\le t\le T$, then
$U_{A}^{\alpha}f$ denotes the $\alpha$-potential of $f\cdot A$,
that is
\[
U_{A}^{\alpha}f(s,x)
=E_{s,x}\int_{s}^{T}e^{-\alpha(t-s)}f(t,X_{t})\,dA_{s,t}\,,\quad
(s,x)\in Q_T.
\]
If $A_{s,t}=t-s$, then we denote $U^{\alpha}_Af$ by $U^{\alpha}f$.

\begin{lm}\label{LU}
For any additive functional $A$ and any $f\in
\mathcal{B}_{b}(Q_{T})$,
\[
U^{\alpha}_{A}(U^{\alpha}f)(s,x)
=E_{s,x}\int_{s}^{T}e^{-\alpha(t-s)}f(t,X_{t})A_{s,t}\,dt
\]
for almost every $(s,x)\in \check{Q}_{T}$.
\end{lm}
\begin{dow}
By the definitions of $\alpha$-potential the fact that $(X,P_{s,x})$ is a Markov process and
Fubini's theorem,
\begin{align*}
U_{K}^{\alpha}(U^{\alpha}f)(s,x)&=E_{s,x}\int_{s}^{T}e^{-\alpha(t-s)}
\left(E_{t,X_{t}}\int_{t}^{T}e^{-\alpha(\theta-s)}
f(\theta,X_{\theta})\,d\theta\right)dA_{s,t}\\
&=E_{s,x}\int_{s}^{T}e^{-\alpha(t-s)}
E_{s,x}\left(\int_{t}^{T}e^{-\alpha(\theta-s)}
f(\theta,X_{\theta})\,d\theta|\mathcal{G}^{s}_{t}\right)dA_{s,t}\\
&=E_{s,x}\int_{s}^{T}\!\!\int_{t}^{T}e^{-\alpha(\theta-s)}
f(\theta,X_{\theta})A_{s,\theta}\,d\theta.
\end{align*}
\end{dow}

\begin{stw}
Let $\mu_{1},\mu_{2}$ be Radon measures such that there exist
additive functionals $K,L$ such that $\mu_{1}\sim K$, $\mu_{2}\sim
L$. If $\mu_{1}\le \mu_{2}$ then $K\le L$ in the sense that
$K_{t^{'},t}\le L_{t^{'},t}$ for every $s\le t^{'}\le t\le T$,
$P_{s,x}$-a.s. for a.e $(s,x)\in\check{Q}_{T}$.
\end{stw}
\begin{dow}
By the assumptions,
\[
U_{K}^{\alpha}f\le U_{L}^{\alpha}f,\quad\alpha\ge0
\]
for every $f\in C^{+}_{0}(\check{Q}_{T})$. Using the theorem on
monotone classes one can show that the above inequalities holds
for any $f\in \mathcal{B}^{+}_{b}(Q_{T})$. In particular, for any
$f\in C^{+}_{0}(\check{Q}_{T})$ and $\alpha\ge0$,
\[
U_{K}^{\alpha}U^{\alpha}f\le U_{L}^{\alpha}U^{\alpha}f.
\]
From this and Lemma \ref{LU} we conclude that for a.e
$(s,x)\in\check{Q}_{T}$,
\[
E_{s,x}f(t,X_{t})K_{s,t}\le E_{s,x}f(t,X_{t})L_{s,t},\quad
t\in[s,T]
\]
for every $f\in C^{+}_{0}(\check{Q}_{T})$. Suppose that $s\le
s^{'}\le t^{'}\le t$. By the above, additivity of $K,L$ and the
Markov property,
\begin{align*}
E_{s,x}f(t',X_{t'})K_{s',t}
&=E_{s,x}f(t',X_{t'})K_{s',t'}+E_{s,x}f(t',X_{t'})K_{t',t}\\
&=E_{s,x}(E_{s',X_{s'}}(f(t',X_{t'})K_{s',t'}))+
E_{s,x}(f(t',X_{t'})E_{t^{'},X_{t'}}K_{t',t})\\
&\le E_{s,x}(E_{s',X_{s'}}(f(t',X_{t'})L_{s',t'}))
+E_{s,x}(f(t',X_{t'})E_{t',X_{t'}}(L_{t',t})])\\
&=E_{s,x}f(t',X_{t'})L_{s',t}.
\end{align*}
By induction, for every $0\le t'\le t_{1}\le\dots\le t_{k}\le t\le
T$ we have
\[
E_{s,x}\prod_{i=1}^{k}f(t_{i},X_{t_{i}})K_{t',t}
=E_{s,x}\prod_{i=1}^{k}f(t_{i},X_{t_{i}})L_{t',t},
\]
from which the lemma follows.
\end{dow}
\begin{wn} If $\mu\sim K$, $\mu\sim L$ then $K=L$.
\end{wn}

It is known (see \cite{Le1,Roz.dec}) that there exist CAF $A$ in
the strict sense and a continuous MAF $M$ in the strict sense such
that
\[
X_t-X_s=M_{s,t}+A_{s,t},\quad t\in[s,T],\quad P_{s,x}\mbox{-}a.s.,
\]
for every $(s,x)\in Q_{\hat T}$, and moreover, $M_{s,\cdot}$ is a
$(\{{\mathcal{G}}^s_t\},P_{s,x})$-square-integrable martingale on
$[s,T]$ with the co-variation  given by
\begin{equation}
\label{eq3.1} \langle M^{i}_{s,\cdot},M^{j}_{s,\cdot}\rangle_t=
\int_s^ta_{ij}(\theta, X_\theta)\,d\theta,\quad t\in[s,T],\quad
i,j=1,...,d,
\end{equation}
while $A_{s,\cdot}$ is a process of $P_{s,x}$-zero-quadratic
variation on $[0,T]$. In particular, $X_{\cdot}-X_s$ is a
$(\{{\mathcal{G}}^s_t\},P_{s,x})$-Dirichlet process in the sense
of F\"ollmer.

Observe that  by (\ref{eq3.1}),
\[
B_{s,t}=\int^t_s\sigma^{-1}(\theta,X_{\theta})\,dM_{s,\theta},\quad
t\in[s,T]
\]
is a $(\{{\mathcal{G}}^s_t\},P_{s,x})$-Wiener process. In
\cite{Le1} it is proved that it has the representation property.
Therefore existence and uniqueness of solutions of (\ref{eq1.6})
follows from known results  for usual BSDEs (see \cite{EKPPQ}),
and moreover, we may  apply Theorem \ref{tw2.1} to RBSDEs with the
Wiener process $B_{s,\cdot}$  defined on the stochastic basis
$(\Omega,\mathcal{G},\{\mathcal{G}^{s}_{t}\}, P_{s,x})$.

We say that a pair $(Y^{s,x},Z^{s,x})$ of
$\{\mathcal{G}^s_t\}$-adapted processes on $[s,T]$ is a solution
of GBSDE (\ref{eq3.04}) if $E_{s,x}\sup_{s\le t\le
T}|Y^{s,x}_{t}|^{2}<\infty$,
$E_{s,x}\int_s^T|Z_t^{s,x}|^2\,dt<\infty$ and (\ref{eq3.04}) is
satisfied $P_{s,x}$-a.s.

Let $S^{s,x}$ be a continuous $\{\mathcal{G}^{s}_{t}\}$ adapted
process. A triple $(Y^{s,x},Z^{s,x},K^{s,x})$ of
$\{\mathcal{G}^s_t\}$-adapted process on  $[s,T]$ is a solution of
RBSDE (\ref{eq1.6}) if $E_{s,x}\sup_{s\le t\le T}
|Y^{s,x}_{t}|^{2}<\infty$,
$E_{s,x}\int_s^T|Z_t^{s,x}|^2\,dt<\infty$,
$E_{s,x}|K^{s,x}_T|^2<\infty$ and (\ref{eq1.6}) is satisfied
$P_{s,x}$-a.s.
\medskip

In the rest of this section we assume that
\begin{enumerate}
\item[(H1)]$\varphi\in\mathbb{L}_{2,\varrho}(\mathbb{R}^{d})$,
\item[(H2)]There exist $M>0$, $\gamma\in\mathbb{L}_{2,\varrho}(Q_{T})$ such
that $|f(t,x,y,z)|\leq |\gamma(t,x)|+M(|y|+|z|)$ for all
$(t,x,y,z)\in[0,T]\times\mathbb{R}^{d}\times\mathbb{R}\times\mathbb{R}^{d}$,
\item[(H3)]$f(t,x,\cdot,\cdot)$ is  continuous for a.e. $(t,x)\in Q_{T}$,
\item[(H4)]$|g(t,x,y)|\le M$
for some $M>0$ and $g(t,x,\cdot)$ is continuous for every
$(t,x)\in Q_{T}$.
\end{enumerate}
\medskip

Let $\mu\in \mathbb{L}_{2}(0,T;H^{-1}_{\varrho})\cap
\mathcal{M}^{+}(Q_{T})$. We say that $u\in \mathcal{W}_{\varrho}$
is a strong solution of the problem (\ref{eqm.0}) if
$u(T)=\varphi$ in $\mathbb{L}_{2,\varrho}(\mathbb{R}^{d})$ and
\[
\langle\frac{\partial u}{\partial t},\eta\rangle_{\varrho,T}
+\langle L_{t}u,\eta\rangle_{\varrho,T} =-\langle f_{u},\eta
\rangle_{2,\varrho,T}-\int_{Q_{T}}\eta g(u)\,d\mu
\]
for every $\eta\in C_{0}^{\infty}(Q_{T})$.
\medskip

Notice that the terminal condition in the above definition is
meaningful since it is known that $\mathcal{W}_{\varrho}\subset
C([0,T],\mathbb{L}_{2,\varrho}(\mathbb{R}^{d}))$.

The following theorem has been proved in \cite{Kl} (see
\cite[Corollary 3.3]{Kl}).
\begin{tw}\label{tw4.6}
Let $h\in\mathcal{W}_{\varrho}$. If $\frac{\partial h}{\partial
t}+L_{t}h=\Phi$ and $\Phi=\alpha_{1}\Phi_1-\alpha_{2}\Phi_2$,
where $\Phi_1,\Phi_2\in\mathbb{L}_2(0,T;H^{-1}_{\varrho}),
\Phi_1,\Phi_2\geqslant0, \alpha_{1},\alpha_{2}\in
\mathcal{B}_{b}(Q_{T}),$ then there exist a cap-quasi continuous
version of $h$, still denoted by $h$, and square-integrable
increasing CAFs $C,R$ such that
\begin{align*}
h(t,X_t)&=h(T,X_T)-\int_{t}^{T}\alpha_{1}(\theta,X_{\theta})\,dC_{s,\theta}
+\int_{t}^{T}\alpha_{2}(\theta,X_{\theta})\,dR_{s,\theta}\\
&\quad-\int_t^T\sigma\nabla
h(\theta,X_\theta)\,dB_{s,\theta},\quad t\in[s,T],\quad
P_{s,x}\mbox{-a.s.}
\end{align*}
for a.e. $(s,x)\in Q_{\hat{T}}$, and if $\mu_1,\mu_2$ are Radon
measures associated with $\Phi_1$ and $\Phi_2$, respectively, then
for a.e. $(s,x)\in Q_{\hat{T}}$,
\begin{equation}
\label{eq3.5} E_{s,x}\int_s^T\xi(\theta,X_\theta)\,dC_{s,\theta}
=\int_s^T\!\!\int_{{\mathbb{R}^d}}\xi(\theta,y)
p(s,x,\theta,y)\,d\mu_1(\theta,y),
\end{equation}
\begin{equation}
\label{eq3.6} E_{s,x}\int_s^T\xi(\theta,X_\theta)\,dR_{s,\theta}
=\int_s^T\!\!\int_{{\mathbb{R}^d}}\xi(\theta,y)
p(s,x,\theta,y)\,d\mu_2(\theta,y)
\end{equation}
for every $\xi \in C_{0}(Q_{T})$ and a.e. $(s,x)\in Q_{T}$.
\end{tw}

The above theorem will be used in the proof of the following
theorem on existence and stochastic representation of strong
solutions of (\ref{eqm.0}) and will play key role in the proof of
Theorem \ref{tw4.1} on existence, approximation and stochastic
representation of strong solutions of the obstacle problem
(\ref{eq1.4}).

\begin{tw}
\label{tw3.3} Assume that (H1)--(H4) are satisfied and
$\mu\in\mathbb{L}_{2}(0,T;H^{-1}_{\varrho})\cap\mathcal{M}^{+}(Q_{T})$.
Then there exists a minimal strong solution $u\in \mathcal{W}_{\varrho}$
of the problem (\ref{eqm.0}). Moreover, the pair
$(u(t,X_{t}),\sigma\nabla u(t,X_{t}))$, $t\in[s,T]$  is a minimal
solution of the GBSDE
\begin{align*}
u(t,X_{t})&=\varphi(X_{T})
+\int_{t}^{T}f_{u}(\theta,X_{\theta})\,d\theta
+\int_{t}^{T}g(u)(\theta,X_{\theta})\,dR_{s,\theta}\\
&\quad -\int_{t}^{T}\sigma\nabla
u(\theta,X_{\theta})\,dB_{s,\theta},\quad t\in[s,T],\quad
P_{s,x}\mbox{-}a.s.,
\end{align*}
where $\mu\sim R.$
\end{tw}
\begin{dow}
First we assume additionally that $f$ is Lipschitz continuous with
respect to $x,y$ uniformly in $t$. Let $(Y^{s,x},Z^{s,x})$ be a
solution of (\ref{eq3.04}).  Existence follows from Proposition
\ref{prog.3} as the assumptions of this theorem are satisfied for
a.e. $(s,x)\in Q_{T}$ (see Corollary 3.5 and  Proposition 4.1 in
\cite{Kl}). Let $g_M(u)=g(u)+M$ so that $g_M(u)\ge 0$, and let
\[
K^{s,x}_{s,t}=\int_{s}^{t}g_M(u)(\theta,X_{\theta},Y^{s,x}_{\theta})\,
dR_{s,\theta},\quad t\in[s,T].
\]
Then we can write (\ref{eq3.04}) in the form
\[
Y^{s,x}_{t}=\varphi(X_{T})
+\int_{t}^{T}f(\theta,X_{\theta},Y^{s,x}_{\theta},Z^{s,x}_{\theta})
+K^{s,x}_T-K^{s,x}_t-MR_{T,t}-\int_{t}^{T}Z^{s,x}_{\theta}\,dB_{s,\theta}.
\]
Let $(Y^{s,x,n},Z^{s,x,n})$ be a solution of the BSDE
\begin{align*}
Y^{s,x,n}_{t}&=\varphi(X_{T})+\int_{t}^{T}
f(\theta,X_{\theta},Y^{s,x,n}_{\theta},Z^{s,x,n}_{\theta})
+K^{s,x,n}_T-K^{s,x,n}_t
-MR_{T,t}\\
&\quad-\int_{t}^{T}Z^{s,x,n}_{\theta}\,dB_{s,\theta},\quad t\in
[s,T],
\end{align*}
where $K^{s,x,n}_{s,t}=\int^{t}_{s}
n(Y^{s,x,n}_{\theta}-Y^{s,x}_{\theta})^{-}\,d\theta$. In much the
same way as in the proof of the approximation result in
\cite[Section 6]{EKPPQ}  (see also \cite{PardouxZhangGRBSDE}) one
can show that
\begin{align}
\label{eq3.08} &E_{s,x}\sup_{s\le t\le T}
|Y^{s,x,n}_{t}-Y^{s,x}_{t}|^{2}
+E_{s,x}\int_{s}^{T}|Z^{s,x,n}_t-Z^{s,x}_t|^{2}\, dt\nonumber\\
&\qquad+E_{s,x}\sup_{s\le t\le T}
|K^{s,x,n}_{s,t}-K^{s,x}_{s,t}|^{2} \rightarrow0
\end{align}
as $n,m \rightarrow \infty$, and
\begin{align}
\label{eq3.09} &E_{s,x}\sup_{s\le t\le T}
|Y^{s,x,n}_{t}|^{2}+E_{s,x}\int_{s}^{T}|Z^{s,x,n}_t|^{2}\,dt
+E_{s,x}|K^{s,x,n}_{s,T}|^{2}\nonumber\\& \qquad \le
CE_{s,x}(|\varphi(X_{T}|^2)+\int_{s}^{T}|\gamma(t,X_t)|^2\,dt
+\sup_{s\le t\le T}|Y^{s,x}_{t}|^{2}+|R_{s,T}|^{2}).
\end{align}
Let us observe now that defining
\begin{equation*}
\tilde{u}(s,x)=E_{s,x}\left(\varphi(X_{T})
+\int_{s}^{T}f(t,X_t,Y^{s,x}_t,Z^{s,x}_t)\,dt
+\int_{s}^{T}g(t,X_t,Y^{s,x}_t)\,dR_{s,t}\right)
\end{equation*}
we get using the Markov property of $\mathbb{X}$ that
\begin{align}
\label{eq4.10}
Y^{s,x}_{t}&=E_{s,x}(Y^{s,x}_t|\mathcal{G}^{s}_{t})\nonumber\\
&=E_{s,x}\left(\varphi(X_{T})+\int_{t}^{T}
f(\theta,X_{\theta},Y^{s,x}_{\theta},Z^{s,x}_{\theta})\, d\theta
+\int_{t}^{T}
g(\theta,X_{\theta},Y^{s,x}_{\theta})\,dR_{s,\theta}
|\mathcal{G}^{s}_{t}\right)\nonumber\\
&=E_{t,X_{t}}\left(\varphi(X_{T})
+\int_{t}^{T}f(\theta,X_{\theta},Y^{s,x}_{\theta},Z^{s,x}_{\theta})\,
d\theta+\int_{t}^{T}
g(\theta,X_{\theta},Y^{s,x}_{\theta})\,dR_{t,\theta}\right)\nonumber\\
&=u(t,X_{t}).
\end{align}
$P_{s,x}$-a.s. for every $t\in [s,T]$. Hence,
\begin{align}
\label{eq4.6} Y^{s,x,n}_{t}&=\varphi(X_{T})+\int_{t}^{T}
f(\theta,X_{\theta},Y^{s,x,n}_{\theta},Z^{s,x,n}_{\theta})\,
d\theta\nonumber\\
&\quad+\int^{T}_{t}n(Y^{s,x,n}_{\theta}
-\tilde{u}(\theta,X_{\theta}))^{-}\,
d\theta-MR_{T,t}-\int_{t}^{T}Z^{s,x,n}_{\theta}\,dB_{s,\theta}
\end{align}
From (\ref{eq4.10}), Corollary 3.5 and  Proposition 4.1 in \cite{Kl} it follows
that $\tilde{u}\in\mathbb{L}_{2,\varrho}(Q_{T})$. In \cite{DPP} it
is proved that that there exists a unique strong solution $u_n$ of
the problem
\begin{equation}\label{eqmm.1}
\frac{\partial u_{n}}{\partial t}+L_{t}u_{n}
=-f_{u_{n}}-n(u_{n}-\tilde{u})^{-}+M\mu,\quad u_{n}(T)=\varphi,
\end{equation}
while from  Theorem \ref{tw4.6} it follows that there is a
cap-quasi continuous version of $u_n$ (still denoted by $u_n$)
such that the pair $(u_n(t,X_t),\sigma\nabla u_n(t,X_t))$,
$t\in[s,T]$, is a solution of (\ref{eq4.6}). Since (\ref{eq4.6})
has a unique solution,
\begin{equation}\label{eqmm.0}
Y^{s,x,n}_{t}=u_{n}(t,X_{t}),\,
t\in[s,T],\,P_{s,x}\mbox{-}a.s.,\quad Z^{s,x,n}_{t}=\sigma\nabla
u_{n}(t,X_{t}),\, \lambda\otimes P_{s,x}\mbox{-}a.s.
\end{equation}

In view of the above we may consider versions of $Y^{s,x,n}$,
$Z^{s,x,n}, Y^{s,x}, Z^{s,x}$ which do not depend on $s,x$. Furthermore, from
(\ref{eq4.10}) it follows that $K^{s,x,n}$, $K^{s,x}$ have
versions not depending on $s,x$. In what follows we consider
versions of the  processes not depending on $s,x$, and
consequently, we drop the superscript $s,x$ in the notation. Write
$d\nu_{n}=n(u_{n}-\tilde{u})^{-}\,dm$. From (\ref{eq3.08}),
(\ref{eq3.09}) and \cite[Proposition 4.1]{Kl} it follows that
\[
\|\nabla u_{n}-\nabla u_{m}\|^{2}_{2\varrho,T} \le
C\int_{Q_{T}}\left(E_{s,x}\int_{s}^{T}
|Z^{n}_{\theta}-Z^{m}_{\theta}|^{2}\,
d\theta\right)\varrho^{2}(x)\,ds\,dx \rightarrow 0
\]
and
\[
\|u_{n}- u_{m}\|^{2}_{2\varrho,T}\le
C\int_{Q_{T}}\left(E_{s,x}\int_{s}^{T}
|Y^{n}_{\theta}-Y^{m}_{\theta}|^{2}\,d\theta\right)\varrho^{2}(x)\,ds\,dx
\rightarrow 0.
\]
Now, if we set $u(t,x)=\lim_{n\rightarrow\infty}u_{n}(t,x)$ if the
limit exists and is finite and $u(s,x)=0$ otherwise, then by the
above, $u_{n}\rightarrow u$ in
$\mathbb{L}_{2}(0,T;H^{1}_{\varrho})$ and
\[E_{s,x}\sup_{s\le t\le T}
|u_{n}(t,X_{t})-u(t,X_{t})|^{2}\rightarrow0
\]
for a.e. $(s,x)\in Q_{\hat{T}}$ which shows that $u$ is cap-quasi
continuous. Now, let $\eta\in C_{0}^{\infty}(Q_{T})$. Then from
the definition of the solution of (\ref{eqmm.1}),

\begin{equation}\label{11}
\langle \frac{\partial u_{n}}{\partial t},\eta\rangle_{\varrho,T}
+\langle L_{t}u_{n},\eta\rangle_{\varrho,T}
=-\langle f_{u_{n}},\eta\rangle_{2,\varrho,T}
-\int_{Q_{T}}\eta\,d\nu_{n}+\int_{Q_{T}}\eta M\,d\mu.
\end{equation}
Hence, by the integration by parts formula,
\begin{align}\label{eqmm.3.5}
\langle u_{n},\frac{\partial\eta}{\partial t}\rangle_{2,\varrho,T}
-\langle L_{t}u_{n},\eta\rangle_{\varrho,T}&=\langle
\varphi,\eta(T)\rangle_{2,\varrho}-\langle
u(0),\eta(0)\rangle_{2,\varrho}+\langle
f_{u_{n}},\eta\rangle_{2,\varrho,T}\nonumber\\
&\quad+\int_{Q_{T}}\eta\,d\nu_{n}-\int_{Q_{T}}\eta M\,d\mu.
\end{align}
By (\ref{eq3.09}), Proposition \ref{stwg.1}, Proposition 4.1 and
Corollary 3.5 in \cite{Kl},
\begin{equation}
\label{eq3.014} \sup_{0\le t\le T}\|u_{n}(t)\|^{2}_{2,\varrho}
+\|\nabla u_{n}\|^{2}_{2,\varrho,T}\le
C(\|\varphi\|^{2}_{2,\varrho} +\|\gamma\|^{2}_{2,\varrho,T}
+\|\mu\|_{\mathbb{L}_{2}(0,T;H^{-1}_{\varrho})}).
\end{equation}
Using this one can check easily that $\{\nu_{n}\}$ is tight.
Therefore without lost of generality we may and will assume that
$\{\nu_{n}\}$ converges weakly to some measure $\nu$.
Consequently, letting $n\rightarrow\infty$ in (\ref{11}) we
conclude that there is functional $\Psi$ on
$C_{0}^{\infty}(\check{Q}_{T})$ such that
\begin{eqnarray}\label{eqmm.5}
\Psi(\eta)+\langle L_{t}u,\eta\rangle_{\varrho,T} =-\langle
f_{u},\eta\rangle_{2,\varrho,T}-\int_{Q_{T}}\eta\,d\nu
+\int_{Q_{T}}\eta M\,d\mu.
\end{eqnarray}
We know that $\nu_{n}\sim K^{n}$, i.e. for a.e. $(s,x)\in
Q_{T}$,
\begin{equation*}
E_{s,x}\int_{s}^{T}\eta(\theta,X_{\theta})\,dK^{n}_{s,\theta}
=\int_s^T\!\!\int_{\mathbb{R}^d}\eta(\theta,y)
p(s,x,\theta,y)\,d\nu_{n}(\theta,y)
\end{equation*}
for every $\eta\in C_{0}(Q_{T})$. Hence, by (\ref{eq3.08}), for
a.e. $(s,x)\in Q_{T}$,
\begin{equation*}
E_{s,x}\int_{s}^{T}\eta(\theta,X_{\theta})\,dK_{s,\theta}
=\int_s^T\!\!\int_{\mathbb{R}^d}\eta(\theta,y)p(s,x,\theta,y)\,d\nu(\theta,y)
\end{equation*}
for every $\eta\in C_{0}(Q_{T})$. Therefore $\nu\sim K$. On the
other hand, by the definition of $K^{n}$,
\begin{equation*}
E_{s,x}\int_{s}^{T}\eta(\theta,X_{\theta})\, dK^{n}_{s,\theta}
=E_{s,x}\int_{s}^{T}\eta(\theta,X_{\theta})
g_M(u_{n})(\theta,X_{\theta})\, dR_{s,\theta}
\end{equation*}
for  $\eta\in C_{0}(Q_{T})$ and a.e. $(s,x)\in Q_{T}$. Using this
and (\ref{eq3.08}), (\ref{eq3.09}) we  conclude that
$g_M(u)\,d\mu\sim K$. By uniqueness, $d\nu=g_M(u)\,d\mu$. Thus,
(\ref{eqmm.5}) takes the form
\[
\Psi(\eta)+\langle L_{t}u,\eta\rangle_{\varrho,T} =-\langle
f_{u},\eta\rangle_{2,\varrho,T}-\int_{Q_{T}}\eta g(u)\,d\mu.
\]
Since $g(u)\,d\mu\in\mathbb{L}_{2}(0,T;H^{-1}_{\varrho})$, using
arguments following (\ref{RM}) we can extend $\Psi$ from
$C_{0}^{\infty}(Q_{T})$ to the functional $\bar{\Psi}$ on
$\mathbb{L}_{2}(0,T;H^{1}_{\varrho})$. Moreover, passing to the
limit in (\ref{eqmm.3.5}) and subtracting (\ref{eqmm.5}) we see
that
\begin{equation*}
\Psi(\eta)=-\langle u,\frac{\partial\eta}{\partial t}
\rangle_{2,\varrho,T}
\end{equation*}
for every $\eta\in C_{0}^{\infty}(\check{Q}_{T})$. Therefore
$\bar\Psi=\frac{\partial u}{\partial t}$ and $u\in
\mathcal{W}_{\varrho}$. Thus, $u$ is a solution of the problem
(\ref{eqm.0}) and, by (\ref{eq3.08})--(\ref{eqmm.0}), the pair
$(u(t,X_{t}),\sigma\nabla u(t,X_{t}))$, $t\in[s,T]$, is a solution
of GBSDE (\ref{eq3.04}). Existence of the minimal solution of
(\ref{eqm.0}) follows now from existence of the minimal solution
of (\ref{eq3.04}) (see Proposition \ref{prog.3}).

We now show how to dispense with the assumption that $f$ is
Lipschitz continuous.  Let $\{f^{n}\}$ be the sequence of
approximations of $f$ considered in Lemma \ref{lmm.1}. Let
$(Y^{s,x,n},Z^{s,x,n})=(Y^n,Z^n)$  be a minimal solution of
(\ref{eq3.04}) with $f$ replaced by $f^n$ and let $u_{n}\in
\mathcal{W}_{\varrho}$ be a solution of the problem
\begin{equation}\label{eqmm.6}
\frac{\partial u_{n}}{\partial t}+L_{t}u_{n}
=-f^{n}_{u_{n}}-g(u_{n})\mu,\quad u_{n}(T)=\varphi.
\end{equation}
From the first part of the proof we know that (\ref{eqmm.0}) is
satisfied. Furthermore, arguing as in the case of usual BSDEs (see
\cite{Lepeltier}) one can show that
\begin{equation}
\label{eq3.16} E_{s,x}\sup_{s\le t\le
T}|Y^{n}_{t}-Y^{m}_{t}|^{2}
+E_{s,x}\int_{s}^{T}|Z^{n}_t-Z^{m}_t|^{2}\,dt \rightarrow 0
\end{equation}
as  $n,m \rightarrow\infty$ and
\begin{align}
\label{eq3.17} &E_{s,x}(\sup_{s\le t\le T}|Y^{n}_{t}|^{2}
+\int_{s}^{T}|Z^{n}_t|^{2}\,dt) \nonumber\\
&\qquad\le CE_{s,x}\left(|\varphi(X_{T})|^2
+\int_{s}^{T}|\gamma(t,X_t)|^2\,dt+|R_{s,T}|^{2}\right).
\end{align}
Set $u(t,x)=\lim u_{n}(t,x)$ if the limit exists and is finite and
$u(s,x)=0$ otherwise. As in in the first part of the proof  we
conclude from (\ref{eq3.16}), (\ref{eq3.17}) and Proposition 4.1
and Corollary 3.5 in \cite{Kl} that $u_{n}\rightarrow u$ in
$\mathbb{L}_{2}(0,T;H^{-1}_{\varrho})$ and $u$ is cap-quasi
continuous. By the definition of the solution of (\ref{eqmm.6}),
\begin{eqnarray}\label{eqmm.7}
\langle \frac{\partial u_{n}}{\partial t},\eta\rangle_{\varrho,T}
+\langle L_{t}u_{n},\eta\rangle_{\varrho,T} =-\langle
f^{n}_{u_{n}},\eta\rangle_{2,\varrho,T} -\int_{Q_{T}}\eta
g(u_{n})\,d\mu
\end{eqnarray}
for every $\eta\in C_{0}^{\infty}(Q_{T})$. The above equality may
be extended to all $\eta\in \mathbb{L}_{2}(0,T;H^{1}_{\varrho})$
(see comments following (\ref{RM})). Moreover,  taking $u_n$ as a
test function  in (\ref{eqmm.7}) and using the properties of the
approximating sequence $\{f^{n}\}$ we conclude that
(\ref{eq3.014}) is satisfied and for every $\eta\in
\mathbb{L}_{2}(0,T;H^{1}_{\varrho})$,
\[
\sup_{n\ge1}|\langle\frac{\partial u_{n}}{\partial t}
,\eta\rangle_{\varrho,T}|<\infty.
\]
This proves that $u\in \mathcal{W}_{\varrho}$ and $\frac{\partial
u_{n}}{\partial t}\rightarrow \frac{\partial u}{\partial t}$
weakly in $\mathcal{W}^{'}_{\varrho}$. By (\ref{eq3.08}),
$u_{n}\rightarrow u$ cap-quasi everywhere. Hence letting
$n\rightarrow\infty$ in (\ref{eqmm.7}) shows that $u$ is a
solution of (\ref{eqm.0}).

Suppose that $v\in \mathcal{W}_{\varrho}$  is another solution of
(\ref{eqm.0}). Then by Theorem \ref{tw4.6} the pair
$(v(t,X_{t}),\sigma\nabla v(t,X_{t}))$, $t\in[s,T]$ is a solution
of GBSDE (\ref{eq3.04}) for a.e. $(s,x)\in Q_{\hat{T}}$. On the
other hand, arguing as in the proof of  Proposition \ref{prog.3}
one can show that $(u(t,X_{t}),\sigma\nabla u(t,X_{t}))$,
$t\in[s,T]$, is the minimal solution of (\ref{eq3.04}) for a.e.
$(s,x)\in Q_{\hat{T}}$. This implies that $u(t,X_{t})\le
v(t,X_{t})$, $t\in [s,T]$ for a.e. $(s,x)\in Q_{\hat{T}}$ which is
equivalent to the fact that $u\le v$ cap-quasi everywhere. Thus,
$u$ is the minimal solution of (\ref{eqm.0}), and the proof is
complete.
\end{dow}

\begin{uw}
{\rm From Theorem \ref{tw3.3} and Remark \ref{rem2.7} it follows
that under assumptions of Theorem \ref{tw3.3} there exists a
maximal solution $u\in \mathcal{W}_{\varrho}$ of (\ref{eqm.0}).  }
\end{uw}

Note that the  stochastic representation of weak solutions of the
problem (\ref{eqm.0}) with $g=0$ was obtained in
\cite{A.Roz.BSDE}.

\nsubsection{Stochastic representation of solutions of the
obstacle problem} \label{sec5}

In this section we consider stochastic homographic approximation
for RBSDEs in a Markovian framework. We assume that the final
condition $\xi$, coefficient $f$ and obstacle $S$ are explicit
functionals of a diffusion associated with the divergence form
operator $L_t$ defined by (\ref{eq1.3}).

We will need the following additional hypotheses.
\begin{enumerate}
\item[(H5)]There is $L>0$ such that
$|f(t,x,y_{1},z_{1})-f(t,x,y_{2},z_{2})|\leq
L(|y_{1}-y_{2}|+|z_{1}-z_{2}|)$ for all
$(t,x)\in[0,T]\times\mathbb{R}^{d}$, $y_{1},y_{2}\in\mathbb{R}$
and $z_{1}, z_{2}\in\mathbb{R}^{d}$,
\item[(H6)]$h\in\mathcal{W}_{\varrho}$, $\varphi(x)\ge h(T,x)$ for a.e.
$x\in\mathbb{R}^d$.
\end{enumerate}

We say that
$\Phi\in{\mathbb{L}}_{2}(0,T;H^{-1}_{\varrho})\cap\mathcal{M}(Q_{T})$
if $\Phi\in{\mathbb{L}}_{2}(0,T;H^{-1}_{\varrho})$ and there
exists $\mu\in\mathcal{M}(Q_{T})$ such that (\ref{RM}) is
satisfied for every $\eta\in C_{0}^{\infty}(Q_{T})$.
\begin{stw}\label{stwo.1}
Let $\mu\in
\mathbb{L}_{2}(0,T;H^{-1}_{\varrho})\cap\mathcal{M}(Q_{T})$ and
let $\mu^{+}-\mu^{-}$ be the Jordan decomposition of $\mu$. Then
$\mu^+,\mu^-\in
\mathbb{L}_{2}(0,T;H^{-1}_{\varrho})\cap\mathcal{M}^{+}(Q_{T}).$
\end{stw}
\begin{dow}
Step 1. First we show that $\mu^{+},\mu^{-}\ll\mbox{Cap}$. Without
lost of generality we can assume that
$\mbox{supp}[\mu]\subset\subset Q_{T}.$ Let $X=\mbox{supp}[\mu]$
and $X=A\cup B$, where $A,B\in \mathcal{B}(Q_{T})$ are from Hahn's
decomposition of signed measure. Let
$A_{\varepsilon},B_{\varepsilon}$ will be compact,
$A_{\varepsilon}\subset A,B_{\varepsilon}\subset B$ and
$|\mu|(A-A_{\varepsilon})<\varepsilon,\,
|\mu|(B-B_{\varepsilon})<\varepsilon.$ Let $K\subset\subset
A_{\varepsilon}$ and $\mbox{Cap}(K)=0$. We will show that
$\mu(K)=0$. Since $\mu\in\mathbb{L}_2(0,T;H^{-1}_{\varrho})$,
$\mu=f-\dyw(\bar{f})$ for some $f,\bar{f}\in
\mathbb{L}_{2,\varrho}(Q_{T})$. Let $A^{\delta}_{\varepsilon}$
will be an open subset of $X$ such that $A_{\varepsilon}\subset
A^{\delta}_{\varepsilon}$ and
$|\mu|(A_{\varepsilon}^{\delta}-A_{\varepsilon})<\delta$. Let
$\xi\in C_{0}^{\infty}(A_{\varepsilon}^{\delta}),\,
\xi_{|A_{\varepsilon}}=1, \xi\ge 0.$ So
$\mu^{\xi}:=\xi\,d\mu=f^{\xi}-\dyw(\bar{f}^{\xi})$ for some
$f^{\xi},\bar{f}^{\xi}\in \mathbb{L}_{2,\varrho}(Q_{T})$. For
$\eta\in C^{\infty}_{0}(Q_{T})$ such that
$\mathbf{1}_{K}\le\eta\le 2$ we have
\begin{align*}
\mu^{\xi}(K)&=\int_{Q_{T}}\mathbf{ 1}_{K}\,d\mu^{\xi}
\le\int_{Q_{T}}\eta\mathbf{ 1}_{A_{\varepsilon}}\,d\mu^{\xi}
=\int_{Q_{T}}\eta\,d\mu^{\xi}+\int_{Q_{T}} (\eta\mathbf{
1}_{A_{\varepsilon}}-\eta)\,d\mu^{\xi}\\
&\le \int_{Q_{T}}\eta\,d\mu^{\xi}+2|\mu|(A^{\delta}_{\varepsilon}
-A_{\varepsilon})\le
C(f^{\xi},\bar{f}^{\xi})\sqrt{\mathcal{E}_{1}(\eta,\eta)}
+2\delta.
\end{align*}
By the above and \cite[Lemma 2.2.7]{Fukushima} we get
\[
\mu(K)=\mu^{\xi}(K)\le C(f^{\xi},\bar{f}^{\xi})Cap(K)
+2\delta=2\delta.
\]
Because $\delta>0$ was arbitrary we get that $\mu^{\xi}(K)=0$ and
hence that $\mu(K)=0$. Similary we obtain that if $K\subset\subset
B_{\varepsilon}$ and $\mbox{Cap}(K)=0$ then $\mu(K)=0$. Let $D\in\mathcal{B}(Q_{T})$ and
$\mbox{Cap}(D)=0$. By \cite[Theorem 2.1.4]{Fukushima} the last
statement is equivalent to $\mbox{Cap}(K)=0$ for every compact
$K\subset D$. Now we have that
\[
\mu(K)=\mu(A\cap K)+\mu(B\cap K)
=\lim_{\varepsilon\rightarrow 0}(\mu(A_{\varepsilon}\cap K)+\mu(B_{\varepsilon}\cap K))=0.
\]
Therefore $\mu(D)=0$. This shows that $\mu \ll\mbox{Cap}$ and as
an immediate consequence that $\mu^{+},\mu^{-}\ll\mbox{Cap}$.

Step 2. Let $\eta\in \mathbb{L}_2(0,T;H^{1}_{\varrho})$. First
assume additionally that $\eta$ is bounded. Since $\mathcal{E}$ is
regular there is a sequence $\{\eta_{n}\}\subset
C_{0}^{\infty}(Q_{T})$ converging to $\eta$ in $\mathcal{E}_{1}$.
By \cite[Theorem 2.1.4]{Fukushima} there exists subsequence
$\{n_{k}\}$ such that $\eta_{n_{k}}\rightarrow \eta$ q.e..  We
know that
\[
\int_{Q_{T}}\eta_{n_{k}}\,d\mu=\langle
f,\eta_{n_{k}}\rangle_{2,\varrho,T}
+\langle\bar{f},\nabla\eta_{n_{k}}\rangle_{2,\varrho,T}
\]
so letting $k\rightarrow\infty$ and using the Lebesgue dominated
convergence theorem we get
\begin{equation}\label{eqo.2}
\int_{Q_{T}}\eta\,d\mu=\langle f,\eta\rangle_{2,\varrho,T}
+\langle\bar{f},\nabla\eta\rangle_{2,\varrho,T}.
\end{equation}
Since for every $\eta\in \mathbb{L}_2(0,T;H^{1}_{\varrho})$ and $c\in\mathbb{R}$ $\eta^{+},\eta^{-},\eta^{+}\wedge c,\eta^{-}\wedge c\in
\mathbb{L}_{2}(0,T;H^{1}_{\varrho})$, using standard arguments we
can show that (\ref{eqo.2}) holds true for any $\eta\in
\mathbb{L}_{2}(0,T;H^{1}_{\varrho})$. In particular, it follows
that $\int_{Q_{T}}\eta\,d|\mu|<\infty$ for
$\eta\in\mathbb{L}_{2}(0,T;H^{1}_{\varrho})$. Since
$\mu=\mu^+-\mu^-\in\mathbb{L}_{2}(0,T;H^{-1}_{\varrho})$ by the
assumption, to prove that $\mu^+,\mu^-\in
\mathbb{L}_{2}(0,T;H^{-1}_{\varrho})$ it suffices to show that
$|\mu|=\mu^++\mu^-\in\mathbb{L}_{2}(0,T;H^{-1}_{\varrho})$, that
is that the functional $|\mu|$ is continuous. Since
$\mathbb{L}_{2}(0,T;H^{1}_{\varrho})$ is a Hilbert space, it
follows from the closed graph theorem that to prove continuity of
$|\mu|$ it  suffices to show that the  $|\mu|$ is closed. But the
last property follows easily from \cite[Theorem 2.1.4]{Fukushima}.
\end{dow}
\medskip

We say that a pair $(u,\mu)$, where $\mu$ is a positive Radon
measure on $Q_T$ and $u\in\mathcal{W}_{\varrho}$, is a strong
solution of the obstacle problem (\ref{eq1.4}) if
\begin{equation}
\label{eq4.1} \langle\frac{\partial u}{\partial
t},\eta\rangle_{\varrho,T}+\langle L_{t}
u,\eta\rangle_{\varrho,T} =\langle
f_u,\eta\rangle_{2,\varrho,T} +\int_{Q_{T}}\eta \,d\mu
\end{equation}
for any $\eta\in\mathcal{F}$,
and
\begin{equation}
\label{eq3.05} u(T)=\varphi,\quad u\geq h\mbox{ on }Q_T,\quad
\int_{Q_T}(u-h)\,d\mu=0.
\end{equation}

It is easily seen that  strong solution of an obstacle problem is
a strong solution of (\ref{eq1.03}) in the variational sense (for
the definition of solution in the variational sense see
\cite{Charrier,Donati,Lions}). Therefore from known results on
uniqueness of variational problems it follows that under (H1),
(H2), (H5), (H6) strong solution of (\ref{eq1.4}) is unique. Let
us observe also that from (\ref{eq4.1}) it follows that
$\mu\in{\mathbb{L}}_{2}(0,T;H^{-1}_{\varrho})$,
which implies that the integral in (\ref{eq3.05}) is well defined.

\begin{tw}\label{tw4.1}
Assume that (H1), (H2), (H5), (H6) are satisfied and
$\frac{\partial h}{\partial t}+L_{t}h\in \mathcal{M}(Q_{T})$. Then
there exists a strong solution $(u,\mu)$ of (\ref{eq1.4}) and if
$u_{n}$, $n\in{\mathbb{N}}$, is a maximal solution of the  Cauchy
problem
\begin{equation}\label{eq.h}
\frac{\partial u_{n}}{\partial
t}+L_{t}u_{n}=-f_{u_{n}}-\mu_n,\quad u_{n}(T)=\varphi
\end{equation}
with
\[
\mu_n=\frac{1}{1+n|u_{n}-h|}\left(\frac{\partial h}{\partial t}
+L_{t}h+f_{h}\right)^{-}
\]
then
\begin{enumerate}
\item[\rm(i)]$u_n\ge h$,  $u_n\searrow u$ a.e.,  $u_n\rightarrow u$
in $\mathbb{L}_2(0,T;H^{1}_{\varrho}) \cap
C([0,T],\mathbb{L}_{2,\varrho}(\mathbb{R}^d))$,
\item[\rm(ii)]$\mu_n\Rightarrow\mu$ and  $\mu_n\rightarrow\mu$ in
$\mathcal{W}^{'}_{\varrho}$.
\end{enumerate}
Moreover for a.e. $(s,x)\in Q_{ \hat{T}}$  there exists a solution
$(Y^{s,x},Z^{s,x},K^{s,x})$ of (\ref{eq1.6}). In fact, for a.e.
$(s,x)\in Q_{\hat{T}}$,
\begin{equation}
\label{eq4.4} Y^{s,x}_t=\tilde u(t,X_t),\quad t\in[s,T],\,\,
P_{s,x}\mbox{-}a.s., \quad Z^{s,x}_t=\sigma\nabla\tilde
u(t,X_t),\quad \lambda\otimes P_{s,x}\mbox{-}a.s.
\end{equation}
for some version $\tilde u$ of $u$ and
\begin{equation}
\label{eq4.5} E_{s,x}\int_s^T\xi(t,X_t)\,dK^{s,x}_t
=\int_s^T\!\!\int_{\mathbb{R}^d}\xi(t,y)p(s,x,t,y)\,d\mu(t,y)
\end{equation}
for every $\xi\in C_{0}(Q_{T})$ and a.e. $(s,x)\in Q_{T}$
\end{tw}
\begin{dow}
Let $\Phi=\frac{\partial h}{\partial t}+L_{t}h+f_{h}$. By the
assumptions on $h$, $\frac{\partial h}{\partial t}+L_{t}h=-
f_{h}+\Phi^{+}-\Phi^{-}$ for some functionals
$\Phi^{+},\Phi^{-}\in{\mathbb{L}}_{2}(0,T;H^{-1}_{\varrho})\cap\mathcal{M}^{+}(Q_{T})$.
Hence, by Theorem \ref{tw4.6}, there exist cap-quasi continuous
versions of $h$ and $u_{n}$ (still denoted by $h,u_{n}$) such that
for a.e. $(s,x)\in Q_{\hat{T}}$,
\begin{align*}
h(t,X_t)&=h(T,X_T)+\int_t^T f_h(\theta,X_\theta)\,d\theta
-C_{t,T}+R_{t,T}\\
&\quad-\int_t^T\sigma\nabla
h(\theta,X_\theta)\,dB_{s,\theta},\quad t\in[s,T],\quad
P_{s,x}\mbox{-}a.s.
\end{align*}
and
\begin{align}
\label{eq4.2}
u_n(t,X_t)&=\varphi(X_T)+\int_t^Tf_{u_n}(\theta,X_\theta)\,d\theta
+\int_t^T\alpha_n(\theta,X_\theta)\,dR_{s,\theta}\\
&\quad-\nonumber \int_t^T\sigma\nabla
u_n(\theta,X_\theta)\,dB_{s,\theta},\quad t\in[s,T],\quad
P_{s,x}\mbox{-}a.s.,
\end{align}
where $\alpha^{n}=\frac{1}{1+n|u_{n}-h|}$ and $C,R$ are  CAFs
associated with $\Phi^{+},\Phi^{-}$(see (\ref{eq3.01})). By Theorem \ref{tw2.1},
$u_n(s,x)=u_n(s,X_s)\ge
u_{n+1}(s,X_s)=u_{n+1}(s,x)$ and $u_n(s,x)=u_n(s,X_s)\ge
h(s,X_s)=h(s,x)$, $P_{s,x}$-a.s. for a.e. $(s,x)\in Q_{\hat{T}}$.
This means that $u_{n}$ is convergent almost everywhere. Set
$u(s,x)=\lim_{n\rightarrow\infty} u_{n}(s,x)$ if the limit exists
and is finite, and $u(s,x)=0$ otherwise. Of course, $u\ge h$ a.e..
By standard calculations, taking $u_{n}$ as a test function in
(\ref{eq.h}), we get
\begin{eqnarray}\label{eq.a}
\|u_n(t)\|^2_{2,\varrho}+\|\nabla u_n\|^2_{2,\varrho,T}\leqslant
C(\|\varphi\|^2_{2,\varrho}+\|g\|^2_{2,\varrho}
+\|\Phi^{-}\|^2_{\mathbb{L}_{2}(0,T;H^{-1}_{\varrho})}).
\end{eqnarray}
By  (\ref{eq2.8}), for a.e. $(s,x)\in Q_{ \hat{T}}$,
\begin{align}
\label{eq5.10} &E_{s,x}|(u_{n}-u_{m})(t,X_{t})|^{2}
+E_{s,x}\int_{s}^{T}|\sigma\nabla(u_{n}-u_{m}|^{2}(\theta,X_{\theta})\,d\theta
\nonumber\\
&\qquad\le E_{s,x}\int_s^T|\alpha^{n}-\alpha^{m}|
|u_{n}-u_{m}|(\theta,X_\theta)\,dR_{s,\theta}
\end{align}
for all $t\in [s,T]$. Multiplying this inequality by $\varrho^{2}$
and using \cite[Proposition 4.1]{Kl} (see also \cite{Bally}) we
obtain
\begin{align*}
&\|u_n(t)-u_m(t)\|^2_{2,\varrho}+\|\sigma\nabla u_n-\sigma\nabla
u_m\|^2_{2,\varrho,T}\\
&\qquad\le CE_{s,\varrho}\int_s^T|(\alpha^{m}-\alpha^{n})|
|u_{n}-u_{m}|(\theta,X_\theta)\,dR_{s,\theta},
\end{align*}
for every $t\in [s,T]$. Here  $E_{s,\varrho}$ denotes the integral
with respect to the measure $P_{s,\varrho}$, where
$P_{s,\varrho}(\cdot)=\int_{\mathbb{R}^d}P_{s,x}(\cdot)\varrho(x)\,dx$.
Next,
\begin{align*}
&E_{s,x}\int_s^T|\alpha^{n}-\alpha^{m}|
|u_{n}-u_{m}|(\theta,X_\theta)\,dR_{s,\theta}\\
&\qquad\le
2E_{s,x}\int_s^T(|h|+|u_{1}|)(\theta,X_\theta)\,dR_{s,\theta}\\
&\qquad \le 2E_{s,x}\sup_{s\le t\le T}
(|h|+|u_{1}|)^{2}(t,X_{t})+2E_{s,x}|R_{s,T}|^{2}.
\end{align*}
By Corollaries 3.4 and 3.6 in \cite{Kl}, $E_{s,\varrho}\sup_{s\le
t\le T}(|h|+|u_1|)^{2}(t,X_t)+E_{s,\varrho}|R_{s,T}|^2<\infty$. On
the other hand, by Theorem \ref{tw2.1}, for a.e. $(s,x)\in
Q_{\hat{T}}$ the right-hand side of (\ref{eq5.10}) tends to zero
when $n,m\rightarrow\infty$. Therefore applying  the Lebesgue
dominated convergence theorem we conclude that $u_{n}\rightarrow
u$ in $\mathbb{L}_2(0,T;H^{1}_{\varrho}) \cap
C([0,T],\mathbb{L}_{2,\varrho}(\mathbb{R}^d))$. Furthermore, if
$\eta\in\mathbb{L}_2(0,T;H^{1}_{\varrho})$, then by (\ref{eq.h}),
\begin{align*}
|\langle \frac{\partial u_{n}}{\partial t},\eta\rangle_{\varrho}|
&\le C\{\|\nabla u_{n}\|_{2,\varrho,T} \|\nabla
\eta\|_{2,\varrho,T} +(\|g\|_{2,\varrho,T}
+\|u_{n}\|_{2,\varrho,T}+\|\nabla u_{n}\|_{2,\varrho,T})
\|\eta\|_{2,\varrho,T}\\
&\quad+\|\Phi^{-}\|_{\mathbb{L}_2(0,T;H^{-1}_{\varrho})}
\|\eta\|_{\mathbb{L}_2(0,T;H^{1}_{\varrho})}\}.
\end{align*}
By (\ref{eq.a}) and the assumptions on the barrier $h$, the
right-hand side of the above inequality is bounded. Hence, by the
Banach-Steinhaus theorem, $\{\frac{\partial u_{n}}{\partial t}\}$
is bounded in $\mathbb{L}_2(0,T;H^{-1}_{\varrho})$.
Consequently, there is a subsequence (still denoted by $\{n\}$)
such that $\frac{\partial u_{n}}{\partial t}\rightarrow
\frac{\partial u}{\partial t}$ weakly in
$\mathbb{L}_2(0,T;H^{-1}_{\varrho})$. By the above
and (\ref{eq.h}) we conclude that
\begin{equation}
\label{eq.he} \frac{\partial u}{\partial
t}+L_{t}u=-f_{u}-\mu,\quad u(T)=\varphi,
\end{equation}
where $\mu$ is a weak limit of $\{\mu_{n}\}$ in
$\mathbb{L}_2(0,T;H^{-1}_{\varrho})$. On the other
hand, passing to the limit in (\ref{eq4.2}) we conclude that there
is an increasing process $K^{s,x}$ on $[s,T]$ such that
$E_{s,x}|K^{s,x}_T|^2<\infty$ and
\begin{align*}
u(t,X_t)&=\varphi(X_T)+\int_t^Tf_{u}(\theta,X_\theta)\,d\theta
+K^{s,x}_T-K^{s,x}_t\\
&\quad-\int_t^T\sigma\nabla
u(\theta,X_\theta)\,dB_{s,\theta},\quad t\in[s,T],\quad
P_{s,x}\mbox{-}a.s.
\end{align*}
for a.e. $(s,x)\in Q_{\hat{T}}$. By the form of above equation we
can drop dependence on $(s,x)$ in notation. This proves
(\ref{eq4.4}) by Theorem \ref{tw2.1}. Formula (\ref{eq4.5})
follows from Theorem \ref{tw4.6}. Now, let us fix
$\eta\in\mathcal{W}_{\varrho}$. From (\ref{eq.h}), (\ref{eq.he})
we get
\begin{align*}
&|\langle \mu_{n}-\mu,\eta \rangle_{\varrho,T}| \le
C\{(\|u_{n}-u\|_{2,\varrho,T} +\|\nabla u_{n}-\nabla u\|)
\|\eta\|_{2,\varrho,T}\\
&\quad+\|\nabla u_{n}-\nabla
u\|_{2,\varrho,T}\|\nabla\eta\|_{2,\varrho,T} +\|\frac{\partial
\eta}{\partial t}
\|_{\mathbb{L}_2(0,T;H^{-1}_{\varrho})}
\|u_{n}-u\|_{\mathbb{L}_2(0,T;H^{1}_{\varrho})}\\&\quad
+\sup_{t\in[0,T]}\|u_{n}(t)-u_{n}(t)\|_{2,\varrho}\|\eta(t)\|_{2,\varrho}\}.
\end{align*}
Since we  know that $u_n\rightarrow u$ in
$\mathbb{L}_2(0,T;H^{1}_{\varrho})
\cap C([0,T],\mathbb{L}_{2,\varrho}(\mathbb{R}^{d}))$,
$\mu_{n}\rightarrow \mu$ in $\mathcal{W}_{\varrho}^{'}$. The proof
is completed by showing that the pair $(u,\mu)$ is a solution of
(\ref{eq1.4}). By what has already been proved all conditions of the definition
of the solution are satisfied but the last condition in (\ref{eq3.05}).
To this end we first note that by the definition
of the solution of RBSDE,
$E_{s,x}\int_{s}^{T}(u-h)(\theta,X_{\theta})\,dK_{s,\theta}=0$ for
a.e. $(s,x)\in Q_{\hat{T}}$. Next, let us observe that formula
(\ref{eq4.5}) can be extended to functions $\xi$ such that one of
the  integrals appearing in (\ref{eq4.5}) is finite. In
particular, by \cite[Corollary 3.6]{Kl}, (\ref{eq4.5}) holds for $\xi=u-h$.
Consequently, for a.e. $s\in [0,T]$ we have
\[
0=\int_{\mathbb{R}^d}\left(E_{s,x}\int_{s}^{T}
(u-h)(\theta,X_{\theta})\,dK_{s,\theta}\right)\,dx\ge
C\int_s^T\!\!\int_{{\mathbb{R}^d}}(u-h)\,d\mu,
\]
the last inequality being a consequence of Aronson's lower
estimate (see \cite{Aron}). This and the fact that any measure in
$\mathbb{L}_2(0,T;H^{-1}_{\varrho})$ vanishes on the
sets of the form $\{t\}\times \mathbb{R}^{d}$ (see comments at the
beginning of  section 4) completes the proof.
\end{dow}
\medskip

In \cite{Palmeri} existence of strong solutions of  variational
problems more general than (\ref{eq1.03}) (with nonlinear operator
$L_t$) is proved by using a different sort  of homographic
approximation. Let us point out that contrary to \cite{Palmeri}
our homographic approximation of $u$ is direct in the sense that
to define the approximating sequence $\{u_n\}$ we need not to
approximate the functional $\left(\frac{\partial h}{\partial t}
+L_{t}h+f_{h}\right)^{-}$ by elements of
$\mathbb{L}_{2,\varrho}(Q_T)$. Secondly, we have proved strong
convergence in $\mathbb{L}_{2,\varrho}(Q_T)$ of gradients  of the
approximating sequence to  the gradient of $u$.

\begin{wn}
\label{wn4.4} Assume that (H1), (H2), (H5), (H6) are satisfied. If
$\frac{\partial h}{\partial t}+L_{t}h\in \mathcal{M}(Q_{T})$,
$(\frac{\partial h}{\partial
t}+L_{t}h-f_{h})^{-}\in\mathbb{L}_{2,\varrho}(Q_{T})$ and
$(u,\mu)$ is a solution of (\ref{eq1.4}), then $d\mu=r\,dm$, where
\begin{equation}
\label{eq3.14} r(t,x)=\alpha(t,x)\left(\frac{\partial h}{\partial
t}+L_{t}h+f_{h}\right)^{-}(t,x)
\end{equation}
for some measurable function $\alpha$  such that
\[
\alpha(t,x)\mathbf{1}_{\{u=h\}}(t,x)=\alpha(t,x),\quad 0\le
\alpha(t,x)\le 1
\]
for a.e. $(t,x)\in Q_T$.
\end{wn}
\begin{dow}
This is an immediate consequence of Theorem \ref{tw4.1} since
$\alpha_{n}$ is bounded in $\mathbb{L}_{2,\varrho}(Q_{T})$ and
$(\frac{\partial h}{\partial t}+L_{t}h+f_{h})^{-} \in
\mathbb{L}_{2,\varrho}(Q_{T})$. The fact that
$\alpha(t,x)\mathbf{1}_{\{u=h\}}(t,x)=\alpha(t,x)$ follows from
uniqueness of the solution of the obstacle problem because the
measure $\mu$ is supported in ${\{u=h\}}$.
\end{dow}

\begin{wn}
Under the assumptions of Corollary \ref{wn4.4} there exists a
version $u$ of the first component of the solution of
(\ref{eq1.4}) such that the triple
\[
\left(u(t,X_{t}),\sigma\nabla u(t,X_{t}),\int_{s}^{t}
r(\theta,X_{\theta})\,d\theta\right),\quad t\in[s,T],
\]
where $r$ is defined by (\ref{eq3.14}), is a solution of  RBSDE
(\ref{eq1.6}) for a.e. $(s,x)\in Q_{\hat{T}}$.
\end{wn}

\begin{uw}
\label{uw5.8} {\rm Let us assume that the operator $L_{t}$ is
symmetric and let assumptions of Theorem \ref{tw4.1} hold. By
Corollary \ref{col2.2}, for every $\eta\in C_{0}^{\infty}(Q_{T})$
we have
\[
0\le E_{s,x}\int_{s}^{T}\eta(t,X_t)\,dK_{s,t} \le
E_{s,x}\int_{s}^{T}\mathbf{1}_{\{u(t,X_t)
=h(t,X_t)\}}\eta(t,X_t)\,dR_{s,t}
\]
for a.e. $(s,x)\in Q_{\check{T}}$. Integrating the above
inequality over $Q_{T}$ with respect to $(s,x)$ and using
(\ref{eq3.6}), (\ref{eq4.5}) and the symmetry of $L_{t}$ we obtain
\[
0\le(\frac{\partial u}{\partial t} +L_{t}u
+f_{u})\le\mathbf{1}_{\{u=h\}}(\frac{\partial h}{\partial
t}+L_{t}h+f_{h})^{-},
\]
i.e. the Lewy-Stampacchia inequality for solutions of the problem
(\ref{eq4.1}), (\ref{eq3.05}) (see \cite{Donati,Palmeri} and
\cite{MM} for the inequality for solutions of elliptic equations).
}
\end{uw}

\end{document}